\documentclass{amsart}

\usepackage{amsmath}
\usepackage{array,amssymb,amsxtra,latexsym,graphics,psfrag,epsfig,ifthen,bbm,bm}
\usepackage[colorlinks,final,backref=page,hyperindex]{hyperref}
\usepackage[all]{xy}
\SelectTips {cm}{}
\usepackage{rotating}
\usepackage[newitem,newenum,neverdecrease]{paralist}
\makeatletter
\if@plflushright
  
\else
  
\fi
\renewcommand{\@asparaenum@}{%
  \expandafter\list\csname label\@enumctr\endcsname{%
    \usecounter{\@enumctr}%
    \labelwidth\z@
    \labelsep.5em
    \leftmargin\z@
    \parsep\parskip
    \itemsep\z@
    \topsep\z@
    \partopsep\parskip
    \itemindent\parindent
    \advance\itemindent\labelsep
    \def\makelabel##1{\upshape ##1}}}
\makeatother
\usepackage{psfrag}
\usepackage{longtable}
\usepackage{xcolor}

\theoremstyle{plain}
\newtheorem{theorem}{Theorem}[section]
\newtheorem{proposition}[theorem]{Proposition}
\newtheorem{lemma}[theorem]{Lemma}

\newtheorem*{theorem*}{Theorem}

\theoremstyle{definition}

\theoremstyle{remark}

\let\map=\xrightarrow
\newcommand{\bdot}{\bm\cdot}

\newcommand{\und}[1]{\underline{#1}}

\newcommand{\id}{\mathrm{id}}

\newcommand{\per}[1]{\begin{pmatrix}#1\end{pmatrix}} 

\newcommand{\Tcp}{\Tc^\bullet}
\newcommand{\Tcc}{\Tc_\bullet}
\newcommand{\Csp}{\Cs^\bullet}
\newcommand{\Csc}{\Cs_\bullet}
\newcommand{\Cspc}{\Cs^\bullet_\bullet}

\newcommand{\Fc}{\mathcal{F}}

\newcommand{\Tc}{\mathcal{T}}

\newcommand{\Mcc}{\mathcal{M}}

\newcommand{\Uc}{\mathcal{U}}

\newcommand{\Nb}{\mathbb{N}}

\let\Kb=\Bbbk

\newcommand{\unit}{I}
\newcommand{\unita}{J}
\newcommand{\unitaa}{K}

\newcommand{\Ss}{\mathsf{Sp}}

\newcommand{\Set}{\mathsf{Set}}
\newcommand{\set}{\mathsf{set}}
\newcommand{\Fset}{\set^{\times}}
\newcommand{\el}{\mathsf{el}}

\newcommand{\Co}[1]{\mathsf{Comon}(#1)}

\newcommand{\Cs}{\mathsf{C}}
\newcommand{\Ds}{\mathsf{D}}
\newcommand{\Is}{\mathsf{I}}

\newcommand{\Ms}{\mathsf{M}}

\newcommand{\rP}{\mathrm{P}}
\newcommand{\rQ}{\mathrm{Q}}
\newcommand{\rB}{\mathrm{B}}

\newcommand{\rL}{\mathrm{L}}
\newcommand{\rU}{\mathrm{U}}
\newcommand{\rI}{\mathrm{I}}
\newcommand{\rJ}{\mathrm{J}}

\newcommand{\qand}{\quad\text{and}\quad}
\newcommand{\qqand}{\qquad\text{and}\qquad}

\begin{document}

\title[Monads on higher monoidal categories]{Monads on higher monoidal categories}

\author[M.~Aguiar]{Marcelo Aguiar}
\address{Department of Mathematics\\
Cornell University\\
Ithaca, NY 14853}
\email{maguiar@math.cornell.edu}
\urladdr{http://www.math.cornell.edu/~maguiar}

\author[M.~Haim]{Mariana Haim}
\address{Centro de Matem\'atica\\
Facultad de Ciencias, Universidad de la Rep\'ublica\\
Montevideo, Uruguay.}
\email{negra@cmat.edu.uy}

\author[I.~L\'opez Franco]{Ignacio L\'opez Franco}
\address{Department of Mathematics and Applications\\
CURE -- Universidad de la Rep\'ublica\\
Tacuaremb\'o s/n, Maldonado, Uruguay.}
\email{ilopez@cure.edu.uy}

\thanks{Aguiar supported in part by NSF grant DMS-140111. Haim and L\'opez~Franco
  supported in part by \textsc{sni--anii} and \textsc{pedeciba}}



\date{\today}

\begin{abstract}
We study the action of monads on categories equipped with several monoidal structures.
We identify the structure and conditions that guarantee that the higher monoidal structure is inherited by the category of algebras over the monad. Monoidal monads and comonoidal monads appear as the base cases in this hierarchy. Monads acting on duoidal categories constitute the next case. We cover the general case of $n$-monoidal categories and discuss several naturally occurring examples in which $n\leq 3$.
\end{abstract}

\maketitle

\setcounter{tocdepth}{1}
\tableofcontents

\section*{Introduction}
\label{s:introduction}

Monoidal categories constitute a natural framework within which certain algebraic structures can be organized, compared and generalized. Classical examples include the category of modules
over a bialgebra (under tensor product over the base field) and the category of modules over a commutative algebra (under tensor product over the algebra). It has long been recognized, on the other hand, that the notion of algebraic theory on a category has a very
simple and powerful generalization in the notion of a monad.
Combining the two settings leads to the consideration of \emph{monoidal monads} and \emph{comonoidal monads}. These monads act on monoidal categories and carry structure that guarantees the existence of a monoidal structure on the category of algebras. The former were considered long ago by Kock \cite{K70,K71a,K71b,K72} and more recently by Seal \cite{S}; the latter in the 90's by McCrudden \cite{McC} and Moerdijk \cite{M}. The monoidal structure on algebras is in one case preserved by the free algebra functor, and by the forgetful functor in the other. Comonoidal monads have allowed in recent years for the generalization of many results pertaining to (bi)algebraic structures \cite{AC,BW,BruVir2007,CLS2010,MesWis2011}.

The notion of \emph{higher monoidal category} specifies suitable compatibility conditions between any given number of monoidal structures on a category.
It originated in \cite{BFSV} and was refined more recently in \cite{AM} and (in the case of two structures) in \cite{G}. The case of \emph{$2$-monoidal} or \emph{duoidal} categories has been the object of much recent study \cite{BM2015,BCZ2013,BGL2014,BS2013,For2004,GL2016,LS2014,St2012}.

This works contributes to this developing theory by studying monads acting on higher monoidal categories. The ultimate goal is to provide conditions that guarantee that the higher monoidal structure is inherited by the category of algebras over the monad. This is reached in Theorem \ref{t:nmonoidal} with the groundwork being laid in earlier sections. Throughout, the notion of multilinearity plays a central role.
The special case in which only two monoidal structures are involved already exhibits most of the interesting phenomena present in the general case. Here the monad may be monoidal with respect to both structures, or comonoidal, or of a mixed nature (monoidal with respect to one and comonoidal with respect to the other). These cases are studied in Sections \ref{s:bimonoidal}, \ref{s:double-mon} and \ref{s:double-comon}. An interesting notion of multilinearity with respect to two monoidal structures is identified in Section \ref{ss:double-multi}.
Sections \ref{s:comonoidal} and \ref{s:monoidal} review the basic cases in which only one monoidal structure intervenes. We devote a certain amount of effort to the case of monoidal monads in Section \ref{s:monoidal}. As mentioned, this case has been dealt with by Seal, building on ealier work of Kock. We add to the discussion by addressing the existence of a multicategory structure
and a colax monoidal structure on the category of algebras over a monoidal monad.
We employ the notion of \emph{representable multicategory} studied by
C.~Hermida~\cite{H}.
The main result (on the monoidal structure of the category of such algebras) is formulated in this context and obtained in Theorem \ref{t:monoidal}. The monoidal structure is defined in terms of certain coequalizers first proposed by Linton \cite{Li}.
An interesting example that combines the constructions of Sections  \ref{s:comonoidal} and \ref{s:monoidal} involves the pointing monad. It leads to the construction of the smash product of bipointed objects in a monoidal category. It is presented in Section \ref{s:pointing} and further developed in Section \ref{ss:pointing-bimon}.
Other examples are presented in the course of the exposition; linear monads are discussed in Sections \ref{ss:linear-monad}, \ref{ss:linear-monoidal-monad}, \ref{ss:linear-bimonoidal-monad}, and \ref{ss:linear-higher}. In Appendix \ref{s:species} we describe certain bimonoidal comonads on the category of species.

\section{Preliminaries}\label{s:preliminaries}

We collect standard material on monads, monoidal and duoidal categories, and multicategories. We refer
to~\cite{BW}, \cite[Chapter 4]{Bor:1994ii}, and \cite[Chapter VI]{McL} for monads and
their algebras, to~\cite{JS}, \cite[Chapters VII and XI]{McL}, and \cite[Chapters 1 and 3]{AM}
for background on monoidal categories, to \cite[Chapter 6]{AM} for duoidal categories, and to~\cite{H} and \cite[Chapters 2 and 3]{Le} for
multicategories and colax monoidal categories.

\subsection{Algebras over a monad}\label{ss:algebra}

Let $(\Tc,\mu,\eta)$ be a \emph{monad} on a category $\Cs$; $\Tc$ is the underlying endofunctor on $\Cs$ and the transformations $\mu$ and $\eta$ are the multiplication and unit of the monad, respectively.

Let $\Cs^{\Tc}$ denote the category of \emph{$\Tc$-algebras} in $\Cs$. The objects are pairs $(V, v)$ where $V$ is an object of $\Cs$ and $v: \Tc V\to V $ is a morphism that is associative and unital. A \emph{morphism of algebras} is required to preserve this structure.

For any $\Tc$-algebra $(V,v)$, the parallel pair of morphisms
\[
\xymatrix@C+5pt{
\Tc^2 V \ar@<0.5ex>[r]^-{\mu_V} \ar@<-0.5ex>[r]_-{\Tc v} & \Tc V
}
\]
is \emph{reflexive}, with common section $\Tc\eta_V$. Moreover, the diagram
\[
\xymatrix@C+5pt{
\Tc^2 V \ar@<0.5ex>[r]^-{\mu_V} \ar@<-0.5ex>[r]_-{\Tc v} & \Tc V\ar[r]^-{v}& V
}
\]
is a \emph{split fork}, and in particular a coequalizer.

A \emph{reflexive coequalizer} is the coequalizer of a reflexive pair. When we speak of coequalizers
\[
\xymatrix@C+5pt{
A \ar@<0.5ex>[r]^-{f} \ar@<-0.5ex>[r]_-{g} & B \ar[r]^-{c}& C,
}
\]
we refer to either the coequalizer object $C$, or the coequalizer map $c:B\to C$, or the whole coequalizer diagram, depending on the context.

\subsection{The free-forgetful adjunction}\label{s:free-forget}

Associated to a monad $\Tc$ on $\Cs$ there is an adjunction
\[
\xymatrix{
\Cs \ar@<0.5ex>[r]^{\Fc^\Tc} & \Cs^\Tc. \ar@<0.5ex>[l]^{\Uc^\Tc}
}
\]
The left adjoint is $\Fc^\Tc X=(\Tc X, \mu_X)$. This is the \emph{free algebra} on the object $X$. The right adjoint forgets the algebra structure: $\Uc^\Tc(V,v)=V$.

Given $f\colon X\to \Uc^\Tc(W,w)$ in $\Cs$, let $\bar f:\Fc^\Tc X\to W$ in $\Cs^\Tc$ be its adjoint, so that $\bar f=w \circ \Tc f$ and $\bar f \circ \eta=f$.

The following result is contained in~\cite[Propositions~4.3.1 and~4.3.2]{Bor:1994ii}.

\begin{lemma}\label{l:creation}
  The functor $\Uc^\Tc$ creates all limits that exist in $\Cs$, and all colimits that exist in $\Cs$ and
  are preserved by both $\Tc$ and $\Tc^2$.
\end{lemma}

The class of reflexive pairs is stable under any functor. Therefore, if $\Tc$ preserves reflexive coequalizers, so does $\Tc^2$.

\begin{lemma}\label{l:creation-reflexive}
Suppose $\Cs$ admits reflexive coequalizers and these are preserved by $\Tc$.
 Then $\Uc^\Tc$ creates reflexive coequalizers.
\end{lemma}

\subsection{Monoidal categories and functors}\label{ss:moncat}

We employ $(\Cs, \otimes, \unit)$ to denote a \emph{monoidal category} and $(\Fc,\varphi,\varphi_0)$ to denote a \emph{monoidal functor} between two monoidal categories $\Cs$ and $\Ds$. The latter involves a natural transformation
\[
\varphi_{X,Y}: \Fc X\otimes \Fc Y \to \Fc(X\otimes Y)
\]
and a map
\[
\varphi_0: I \to \Fc I
\]
from the unit object of $\Ds$ to the image of the unit object of $\Cs$ under $\Fc$, satisfying certain axioms. A \emph{comonoidal functor} $(\Fc,\psi,\psi_0)$ is the dual notion. When $(\Fc,\varphi,\varphi_0)$ is monoidal and $\varphi$ and $\varphi_0$ are invertible, $(\Fc,\varphi^{-1},\varphi_0^{-1})$ is comonoidal, and vice versa. In this case we say $\Fc$ is \emph{strong}.

We often treat monoidal categories as if they were \emph{strict}.

Iterating the structure $(\varphi,\varphi_0)$ of a monoidal functor $\Fc$ one obtains maps
\[
\varphi_{X_1,\dots,X_n}: \Fc X_1 \otimes \dots \otimes \Fc X_n \to \Fc(X_1\otimes \dots \otimes X_n)
\]
for any sequence of objects $X_1,\dots,X_n$. When $n=1$, this map is the identity. When $n=0$, this is $\varphi_0$.

\subsection{Duoidal categories}\label{ss:duoidal}

Let $(\diamond, \unit)$ and $(\star, \unita)$ be two monoidal structures on the same category $\Cs$.  Suppose there is in addition a natural transformation (the \emph{interchange law})
\[
\zeta_{A,B,C,D} :
(A \star B) \diamond (C \star D) \to (A \diamond C) \star (B \diamond D)
\]
and three morphisms (the \emph{unit maps})
\[
 \unit \to \unit \star \unit \qquad
  \unita \diamond \unita \to \unita \qquad \unit \to \unita
\]
satisfying the axioms given in \cite[Definition 6.1]{AM}. We state one of them:
\begin{equation}\label{eq:iass}
\begin{gathered}
\xymatrix@C-6pc@R+1pc{
& (A \star X) \diamond (B \star Y) \diamond (C \star Z) \ar[dl]_{\zeta \diamond \id} \ar[dr]^{\id \diamond \zeta} & \\
\bigl((A \diamond B) \star (X \diamond Y)\bigr) \diamond (C \star Z) \ar[dr]_{\zeta} & &
(A \star X) \diamond \bigl((B \diamond C) \star (Y \diamond Z)\bigr) \ar[dl]^{\zeta}\\
& (A \diamond B \diamond C) \star (X \diamond Y \diamond Z)  &
}
\end{gathered}
\end{equation}
Endowed with such structure, $\Cs$ is a \emph{$2$-monoidal category}. The term \emph{duoidal category} is also employed.

The interchange law and the unit maps are not required to be invertible. We indicate the direction of these maps by listing the structures in order and saying that $(\Cs,\diamond,\star)$ is a duoidal category.

Suppose a monoidal category $(\Cs,\otimes,\unit)$ carries a braiding $\beta$.
In this case, $\Cs$ is duoidal with both structures $(\diamond,\unit)$ and $(\star,\unita)$ equal to $(\otimes,\unit)$ and with interchange law
\[
\zeta_{A,B,C,D}=\id_A \otimes \beta_{B,C} \otimes \id_D.
\]
The unit morphisms are canonical isomorphisms. In this case, all structure maps are invertible.
A version of the Eckmann-Hilton argument due to Joyal and Street \cite[Proposition 6.11]{AM} states the converse: if all structure maps are invertible, the duoidal category must arise from a braided monoidal category in this manner.

For a simple example of a different nature, let $G$ be a monoid (in the cartesian category of sets). The category of $G$-graded sets is duoidal. The objects are sequences $S=(S_g)_{g\in G}$ of sets. The monoidal structures $S\diamond T$ and $S\star T$ are given by
\[
(S\diamond T)_g = \coprod_{xy=g} S_x\times T_y \qand (S\star T)_g=S_g\times T_g.
\]
The interchange law is an inclusion.

Additional examples of duoidal categories are given in \cite[Section 6.4]{AM} and later in this paper.

\subsection{Colax monoidal categories}\label{ss:colax}

 A \emph{colax monoidal} structure on a category $\Cs$ consists of a sequence of functors
 \[
 \otimes^n:\Cs^n\to\Cs,\ n\geq 0,
 \]
a natural transformation $\epsilon:\otimes^1\to\id$,
and a family of natural transformations $\alpha^{m_1,\dots,m_n}$, $m_1,\dots,m_n\geq 0$, as follows. For each $i=1,\ldots,n$, let
\[
\und{X}^i = (X_1^i,\dots,X_{m_i}^i)
\]
be a sequence of length $m_i$ of objects of $\Cs$. Let
\[
\und{X} = (X_1^1,\dots,X_{m_n}^n)
\]
denote the concatenation of the $n$ sequences $\und{X}^1,\dots,\und{X}^n$. The transformations map
\begin{equation}\label{eq:colax}
\alpha^{m_1,\dots,m_n}:  \otimes^{m_1+\dots+m_n}(\und{X}) \to
  \otimes^n\bigl(\otimes^{m_1}(\und{X}^1),
  \dots, \otimes^{m_n}(\und{X}^n)\bigr).
\end{equation}
They are subject to certain axioms: see \cite[Definition 3.1.1]{Le} (for the dual axioms).

When the transformation $\epsilon$ is invertible we say that the structure is \emph{normal}.
We denote such structure by $(\Cs, \und{\otimes},\und{\alpha})$, where $\und{\otimes}$ stands for the sequence
$\{\otimes^n\}_{n\geq 0}$ and $\und{\alpha}$ for the family $\{\alpha^{m_1,\dots,m_n}\}_{m_i\geq 0}$. Note that $\otimes^0$ amounts to the choice of an object of
$\Cs$.

Every monoidal structure $(\otimes, \unit)$ on a category $\Cs$ gives rise to a normal colax monoidal structure on the same category with
\[
  \otimes^n(X_1,\dots,X_n)=X_1\otimes \cdots \otimes X_n,
  \quad
  \otimes^0=I.
\]
The transformations \eqref{eq:colax} are invertible and constructed from the
associativity and unit constraints of $\Cs$. The transformation $\epsilon$ is the identity.

\subsection{Multicategories}\label{s:multicat}

A \emph{multicategory} $\Ms$ possesses objects and \emph{multimaps}. The domain of a multimap is a finite sequence of objects and the codomain is a single object.
The notation $f: (X_1,\dots, X_n) \to Y$ indicates that $f$ is multimap from $(X_1, \dots, X_n)$ to $Y$ in $\Ms$. In this case, we say that $f$ is an $n$-map.
Multimaps can be suitably composed (in a tree-like manner).
Each object has an associated identity (a $1$-map with domain and codomain equal to the given object).

A multicategory $\Ms$ has an \emph{underlying category} $\Uc(\Ms)$. It has the same objects as $\Ms$ and the morphisms are the $1$-maps of $\Ms$.

A \emph{universal multimap} for the sequence $(X_1, \dots X_n)$ is a multimap
\[
u: (X_1,\dots,X_n) \to U
\]
in $\Ms$ through which every other multimap of the same domain factors. In other words, for every multimap $f:(X_1,\dots,X_n) \to Y$ in $\Ms$, there exists a unique $1$-map $\hat f:U\to Y$ making the diagram
\[
\xymatrix@C+10pt{
(X_1,\dots ,X_n) \ar[r]^-f \ar[d]_{u} &Y \\
U \ar[ru]_-{\hat f}
}
\]
commutative. For a given sequence of objects, universal multimaps may or may not exist. When they do, they are unique up to isomorphism. For any object $X$, $\id_X$ is a universal $1$-map for $X$.

\subsection{Representable multicategories}\label{ss:multicat-colax}

A colax monoidal category $(\Cs,\und{\otimes},\und{\alpha})$ gives rise to a multicategory $\Mcc(\Cs)$.  The objects of $\Mcc(\Cs)$ are those of $\Cs$.  A multimap $X_1,\dots, X_n \to Y$ is a morphism in $\Cs$ of the form
\[
f: \otimes^n(X_1,\dots,X_n) \to Y.
\]
Composition of multimaps employs the transformations $\und{\alpha}$. The identity of $X$ in $\Mcc(\Cs)$ is $\epsilon_X:\otimes^1 X\to X$.

There is a canonical functor $\Cs\to \Uc\Mcc(\Cs)$ which sends
\[
X\map{f} Y \text{ \ to \ } \otimes^1 X\map{\epsilon_X} X \map{f} Y.
\]
When $\Cs$ is normal, in fact $\Cs\cong \Uc\Mcc(\Cs)$.

A multicategory $\Ms$ is (colax) \emph{representable} if
there exists a normal colax monoidal structure on $\Uc(\Ms)$ such that
$\Ms\cong \Mcc\Uc(\Ms)$ under an isomorphism of multicategories
that restricts to the canonical isomorphism between the underlying categories.

The following lemma is a variant of the results of~\cite[Section 9]{H} on representable multicategories.

\begin{lemma}\label{l:inducedmulticat}
A multicategory $\Ms$ is representable if and only if $\Ms$ admits universal multimaps. In this case, the colax monoidal structure on $\Uc(\Ms)$ is monoidal if and only if the class of universal multimaps is closed under composition.
\end{lemma}

The codomain of the universal multimap $u$ corresponding to
  the tuple $X_1,\dots,X_n$ is $\otimes^n(X_1, X_2, \dots , X_n)$. The transformations $\epsilon$ and $\und{\alpha}$ are defined by universality. In particular,
\[
\xymatrix@R-5pt{
X \ar[d] \ar[r]^-{\id_X} & X\\
\otimes^1 X \ar[ru]_{\epsilon_X}
}
\]
and since $\id_X$ is itself universal, we have that $\epsilon_X$ is invertible, so the colax structure is normal.

\section{Comonoidal monads}\label{s:comonoidal}
We review the definition of comonoidal monad and its basic properties. References include~\cite{McC,M,Szl}. (The term \emph{opmonoidal monad} is used in \cite{McC,Szl} for these objects, and the term \emph{Hopf monad} is used in \cite{M}.)

\subsection{Comonoidal monads and their algebras}\label{ss:algcomonoidal}

Let $(\Cs, \otimes, \unit)$ be a monoidal category. A monad $(\Tc, \mu, \eta)$ on $\Cs$
 is \emph{comonoidal} if the functor $\Tc$ is equipped with a comonoidal
  structure $(\psi,\psi_0)$ for which $\mu$ and $\eta$ are morphisms
 of comonoidal functors. In other words, a comonoidal monad is a monad in the $2$-category of monoidal categories and comonoidal functors. 


Let $\Tc$ be such a monad. Given two $\Tc$-algebras $(V,v)$ and $(W,w)$, we may turn $V\otimes W$ into an algebra by means of
\[
\Tc(V\otimes W) \map{\psi_{V,W}} \Tc V \otimes \Tc W \map{v\otimes w} V\otimes W.
\]
The unit object of $\Cs$ is also an algebra by means of
\[
\Tc\unit \map{\psi_0} \unit.
\]

\begin{theorem}\label{t:comonoidal}
 The preceding defines a monoidal structure on the
 category $\Cs^\Tc$ of $\Tc$-algebras.  It is such that the forgetful functor
 $U^\Tc: \Cs^\Tc \to \Cs$ is strict monoidal.
\end{theorem}

This is  \cite[Proposition 1.1]{McC} and \cite[Theorem 7.1]{M}.
 
The compatibility between the monad and comonoidal structures on $\Tc$ expresses precisely that $\psi_0$ and $\psi_{X,Y}$ are morphisms of $\Tc$-algebras, for any objects $X$ and $Y$ of $\Cs$.

Theorem \ref{t:monoidal} admits the following converse. Suppose $\Tc$ is a monad and that $\Cs^\Tc$ carries a monoidal structure for which $\Uc^\Tc$ is strong monoidal. Then $\Fc^\Tc$ is comonoidal by means of the mate structure \cite{Kel74} (see, for example, \cite[Section 3.9]{AM}). Then $\Tc=\Fc^\Tc\Uc^\Tc$ is comonoidal and the structure is such that $\psi_0$ and $\psi_{X,Y}$ are morphisms of $\Tc$-algebras. Hence $\Tc$ is a comonoidal monad.

\subsection{Example: linear monads}\label{ss:linear-monad}

Let $H$ be a bimonoid in a $2$-monoidal category $(\Cs,\diamond,\star)$. In particular, $H$ is a monoid in $(\Cs,\diamond)$ and a comonoid in $(\Cs,\star)$. Consider the endofunctor $\Tc$ on $\Cs$ defined by
\[
\Tc X = H\diamond X.
\]
The monoid structure of $H$ results in a monad structure on $\Tc$. Monads arising from monoids in this manner are called \emph{linear}. The comonoid structure of $H$, on the other hand, may be employed to turn $\Tc$ into a comonoidal monad on the monoidal category $(\Cs,\star)$. The comonoidal structure $\psi_{X,Y}$ is
\[
\Tc(X\star Y) = H\diamond (X\star Y) \map{\Delta\diamond\id} (H\star H)\diamond (X\star Y) \map{\zeta}
(H\diamond X)\star(H\diamond Y) = \Tc X \star \Tc Y.
\]
The map $\psi_0$ is defined similarly, in terms of the counit of $H$ and one of the unit maps of $\Cs$. The bimonoid axioms imply those for a comonoidal monad.

In this case, $\Tc$-algebras are modules over the monoid $H$ in $(\Cs,\diamond)$. Theorem \ref{t:comonoidal} affords a monoidal structure on this category. This recovers
 \cite[Proposition~6.39]{AM}. The special case when $H$ is a bimonoid in a braided monoidal category is classical.

\subsection{Example: monads on cartesian categories}\label{ss:cartesian}

Let $\Cs$ be a category admitting all finite products. We may view it as a monoidal category in which the tensor product is the categorical product and the unit object is a terminal object.
Categories of this form are called \emph{cartesian}.

Any functor between cartesian categories carries a unique comonoidal structure. (See, for instance, \cite[Example 3.19]{AM}.) In addition, this structure is preserved by any natural transformation between two such functors. It follows that any monad on a cartesian category admits a unique comonoidal structure.

Theorem \ref{t:comonoidal} says that products of algebras over such a monad are given by products in $\Cs$. This well-known fact also follows from Lemma \ref{l:creation}.





\section{Monoidal monads}\label{s:monoidal}

We review the definition of monoidal monad and go over the construction of a monoidal structure on the category of algebras in some detail. We discuss the connection between Linton's coequalizers, which define this structure, and universal multilinear maps (multimaps in a multicategory whose objects are algebras over the monoidal monad). The multicategory is defined under no further assumptions. When universal multilinear maps exist, the multicategory is representable and the category of algebras carries a colax monoidal structure. Under additional assumptions, this structure is in fact monoidal.

\subsection{Monoidal monads and Linton's coequalizers}\label{ss:monad-Linton}

Let $(\Cs, \otimes, \unit)$ be a monoidal category. A monad $(\Tc, \mu, \eta)$ on $\Cs$
 is \emph{monoidal} if the functor $\Tc$ is equipped with a monoidal
 structure $(\varphi,\varphi_0)$ for which $\mu$ and $\eta$ are morphisms
 of monoidal functors. In other words, a monoidal monad is a monad in the $2$-category of monoidal categories and monoidal functors. 

Let $\Tc$ be such a monad. Given $\Tc$-algebras $(V_i,v_i)$, $i=1,\ldots,n$, consider the following parallel pair of morphisms:
\begin{equation}\label{eq:pair}
\xymatrix{
	\Tc(\Tc V_1 \otimes \dots \otimes \Tc V_n) \ar[r]^-{\Tc\varphi}
	\ar@/_1pc/[rr]_-{\Tc(v_1\otimes \dots \otimes v_n)} &
	\Tc^2(V_1\otimes  \dots \otimes V_n) \ar[r]^-{\mu}
	& \Tc(V_1\otimes \dots \otimes V_n).
}
\end{equation}
It is reflexive: the map
\[
\xymatrix@C+60pt{
\Tc(\Tc V_1 \otimes \dots \otimes \Tc V_n) &
 \Tc(V_1\otimes \dots \otimes V_n)
 \ar[l]_-{\Tc(\eta_{V_1} \otimes \dots \otimes \eta_{V_n})}
}
\]
is a common section. Note that each of these objects is a free $\Tc$-algebra and each of these maps is a morphism of $\Tc$-algebras.

We refer to an equalizer of a pair of the form \eqref{eq:pair} as a \emph{Linton coequalizer}.
They were considered in \cite{Li}. We will consider such equalizers both in the category $\Cs$ and in the category $\Cs^\Tc$ of $\Tc$-algebras.

If all Linton coequalizers exist in $\Cs$ and they are preserved by $\Tc$ and $\Tc^2$, then they also exist in $\Cs^\Tc$ and they are calculated as in $\Cs$, according to Lemma \ref{l:creation}.
The same conclusion holds if $\Cs$ admits all reflexive coequalizers and these are preserved by $\Tc$, in view of Lemma \ref{l:creation-reflexive}.

In the case of free algebras, Linton coequalizers always exist, according to the following result.
For $i=1,\dots,n$, let $X_i$ be an object of $\Cs$.

\begin{lemma}\label{l:linton-free}
The diagram
\[
\xymatrix@C+10pt{
\Tc(\Tc^2 X_1\otimes\dots\otimes\Tc^2 X_n)
\ar@<0.5ex>[r]^-{\mu\circ\Tc\varphi}  \ar@<-0.5ex>[r]_-{\Tc(\mu\otimes\dots\otimes\mu)} &
\Tc(\Tc X_1\otimes\dots\otimes\Tc X_n)
\ar[r]^-{\mu\circ\Tc\varphi} &
\Tc(X_1\otimes\dots\otimes X_n)
}
\]
is a split cofork (and hence a coequalizer) in $\Cs^\Tc$.
\end{lemma}
\begin{proof}
The splitting data is
\[
\xymatrix@C+15pt{
\Tc(\Tc^2 X_1\otimes\!\cdots\!\otimes\Tc^2 X_n)  &
\Tc(\Tc X_1\otimes\!\cdots\!\otimes\Tc X_n)
\ar[l]_-{\Tc(\Tc\eta\otimes\dots\otimes\Tc\eta)} &
\Tc(X_1\otimes\!\cdots\!\otimes X_n). \ar[l]_-{\Tc(\eta\otimes\dots\otimes\eta)}
}
\]
\end{proof}


\subsection{Multilinearity}\label{ss:multilinear}

Let $(\Tc, \varphi)$ be a monoidal monad on $\Cs$.
For $i=1,\dots,n$, let $(V_i,v_i)$ be a $\Tc$-algebra, as before, and let $(W,w)$ be another such algebra. A morphism $f:V_1\otimes \cdots \otimes V_n \to W$ in $\Cs$ is said to be \emph{$n$-linear} (with respect to $(\Tc, \varphi)$) if the following diagram commutes.
\begin{equation}\label{eq:multilinear}
\begin{gathered}
\xymatrix{
\Tc V_1 \otimes \cdots \otimes \Tc V_n \ar[r]^-{\varphi} \ar[dd]_-{v_1 \otimes \cdots \otimes v_n}& \Tc(V_1\otimes \cdots \otimes V_n) \ar[d]^-{\Tc f}  \\
& \Tc W \ar[d]^-w\\
V_1\otimes \cdots \otimes V_n \ar[r]_-f & W
}
\end{gathered}
\end{equation}
We refer to such maps $f$ collectively as \emph{$\Tc$-multilinear}. When $n=2$, we
say $f$ is  $\Tc$-bilinear, and when $n=1$ we say it is $\Tc$-linear. A $\Tc$-linear map is the same as a morphism of $\Tc$-algebras. A $0$-linear map is the same as a morphism $\unit\to W$ in $\Cs$: the commutativity of \eqref{eq:multilinear} is in this case automatic, it follows from the fact that $\varphi_0=\eta_\unit$ which is one of the conditions for a monoidal monad.

Define a multicategory $\Cs^{\Tc,\varphi}$ as follows. The objects are $\Tc$-algebras.
A multimap $V_1,\cdots,V_n \to W$ in $\Cs^{\Tc,\varphi}$ is a $\Tc$-multilinear map $f:V_1\otimes \cdots \otimes V_n \to W$ in $\Cs$. Such maps are particular multimaps in the multicategory $\Mcc(\Cs,\otimes)$ associated to the monoidal category $(\Cs,\otimes)$ (Section \ref{ss:multicat-colax}).

\begin{lemma}\label{l:multisubcat}
The forgetful functor $\Uc^\Tc:\Cs^\Tc\to\Cs$ extends to a faithful morphism of multicategories
\[
\Cs^{\Tc,\varphi}\to\Mcc(\Cs,\otimes).
\]
\end{lemma}
\begin{proof}
Every identity is linear. We omit the verification of the fact that  multilinear maps are closed under composition in $\Mcc(\Cs,\otimes)$.
\end{proof}

Since linear maps are precisely algebra morphisms, we have $\Uc(\Cs^{\Tc,\varphi})=\Cs^\Tc$.

\subsection{Kock's criterion}\label{ss:kock}

The notion of multilinearity originates in work of Linton \cite{Li} and Kock \cite{K71a}. The following criterion for bilinearity is given in \cite[Theorem 1.1]{K71a}. There is an evident extension to multilinearity. 

\begin{lemma}\label{l:kock}
Let $A$, $B$, $C$ be $\Tc$-algebras. A map $f:A\otimes B\to C$ in $\Cs$ is $\Tc$-bilinear if and only if the following diagrams commute.
\[
\xymatrix@-10pt{
\Tc A \otimes B \ar[r]^-{\id\otimes\eta}  \ar[dd]_-{a \otimes \id}& \Tc A\otimes\Tc B \ar[r]^-{\varphi} &\Tc(A\otimes B) \ar[d]^-{\Tc f}  \\
& & \Tc C \ar[d]^-c\\
A\otimes B \ar[rr]_-f & & C
}
\qquad
\xymatrix@-10pt{
A \otimes \Tc B \ar[r]^-{\eta\otimes\id} \ar[dd]_-{\id \otimes b}& \Tc A \otimes \Tc B \ar[r]^-{\varphi} & \Tc(A\otimes  B) \ar[d]^-{\Tc f}  \\
& & \Tc C \ar[d]^-c\\
A\otimes B \ar[rr]_-f & & C
}
\]
\end{lemma}

We employ Kock's criterion in a few examples later.

\subsection{Linton's coequalizers and universal multilinear maps.}\label{ss:Linton-univ-multi}

We continue to work with a monoidal monad $(\Tc,\varphi)$ on the monoidal category
  $(\Cs, \otimes)$. Let $(V_i,v_i)$, $i=1,\dots,n$, and $(W,w)$ be $\Tc$-algebras.
  Given $f:V_1\otimes \dots \otimes V_n \to W$ in $\Cs$, let $\bar{f}$ in $\Cs^\Tc$ be its adjoint:
\[
\bar f:\Tc(V_1\otimes\dots\otimes V_n) \map{\Tc f} \Tc W \map{w} W.
\]
The map $f$ is recovered as
\[
f: V_1\otimes \dots \otimes V_n \map{\eta} \Tc(V_1\otimes\dots\otimes V_n) \map{\bar{f}} W.
\]
We first look at the behavior of this correspondence under post-composition.
Suppose $h:(W,w)\to(R,r)$ is a morphism of algebras.

\begin{lemma}\label{l:post-tilde}
We have
\[
h\circ\bar{f} = \overline{h\circ f}.
\]
\end{lemma}
\begin{proof}
Precomposing either side with $\eta$ yields $h\circ f$.
\end{proof}

Multilinearity of $f$ and $\bar{f}$ being a Linton coequalizer are closely related properties.

\begin{lemma}\label{l:coeqmult}
 We have that
  \begin{enumerate}[(i)]
  \item\label{item:5}
  $f$ is $\Tc$-multilinear if and only if $\bar{f}$ coequalizes the pair~\eqref{eq:pair};
  \item\label{item:6}
  $f$ is a universal multimap in $\Cs^{\Tc,\varphi}$ if and only if $\bar{f}$ is the coequalizer of the pair~\eqref{eq:pair} in $\Cs^\Tc$.
  \end{enumerate}
\end{lemma}
\begin{proof}
Multilinearity for $f$ is the equality
\[
w\circ \Tc f \circ \varphi = f \circ (v_1\otimes \dots \otimes v_n),
\]
while the coequalizing condition for $\bar f$ is
\[
w\circ \Tc f \circ \mu \circ \Tc\varphi =w\circ \Tc f \circ \Tc(v_1\otimes \dots \otimes v_n).
\]
The former equality is equivalent (by adjointness) to
\[
w\circ \Tc(w\circ \Tc f \circ \varphi)=w\circ \Tc\left (f \circ (v_1\otimes \dots \otimes v_n)\right),
\]
so to prove \eqref{item:5} it suffices to verify that
\[
w\circ \Tc(w\circ \Tc f \circ \varphi)=w \circ \Tc f \circ \mu \circ \Tc \varphi .
\]
This holds by functoriality of $\Tc$, naturality of $\mu$, and associativity of $w$.

To prove \eqref{item:6} we may assume that $f$ is $\Tc$-multilinear (this is either an assumption, or a consequence, in view of \eqref{item:5}). Part \eqref{item:5} tells us that then $\bar{f}$ coequalizes the pair~\eqref{eq:pair}.

Let $(R,r)$ be a $\Tc$-algebra. Let $\Cs^{\mathrm{coeq}}\bigl((V_1,v_1),\dots,(V_n,v_n);(R,r)\bigr)$ denote the set of algebra morphisms
\[
\Tc(V_1\otimes\dots\otimes V_n) \to R
\]
that coequalize the pair~\eqref{eq:pair}. By \eqref{item:5} applied to $(R,r)$, we have a bijection
\[
\Cs^{\Tc,\varphi}\bigl((V_1,v_1),\dots,(V_n,v_n);(R,r)\bigr)\to
\Cs^{\mathrm{coeq}}\bigl((V_1,v_1),\dots,(V_n,v_n);(R,r)\bigr),  \quad g\mapsto \bar{g}.
\]
Consider the following diagram
\[
\xymatrix@C-50pt{
\Cs^{\Tc,\varphi}\bigl((V_1,v_1),\dots,(V_n,v_n);(R,r)\bigr) \ar[rr]^-{\cong} & & \Cs^{\mathrm{coeq}}\bigl((V_1,v_1),\dots,(V_n,v_n);(R,r)\bigr)\\
& \Cs^{\Tc}\bigl((W,w),(R,r)\bigr) \ar[lu]^-{(-)\circ f} \ar[ru]_-{(-)\circ\bar{f}} &
}
\]
in which the diagonal arrows are given by precomposition with $f$ and $\bar{f}$, respectively.
Lemma \ref{l:post-tilde} says precisely that the diagram commutes. Now, a multilinear map $f$ as in the statement is universal if and only if $(-)\circ f$ is a bijection
for all algebras $(R,r)$. Since the horizontal map is a bijection, this is equivalent to $(-)\circ\bar{f}$ being a bijection for all $(R,r)$. Given that $\bar{f}$ coequalizes, this is in turn equivalent to $\bar{f}$ being the coequalizer of the pair~\eqref{eq:pair}.
\end{proof}

\subsection{Free algebras and multilinearity}\label{ss:free-multi}

We turn to the construction of a morphism of multicategories $\Mcc(\Cs,\otimes)\to \Cs^{\Tc,\varphi}$ which agrees with the free algebra functor on the underlying categories.

\begin{lemma}\label{l:phi-multi}
Let $n\geq 0$ and $X_i$ be an object of $\Cs$, $i=1,\dots,n$.
The structure map
\[
\varphi_{X_1,\dots,X_n}: \Tc X_1\otimes\dots\otimes\Tc X_n \to \Tc(X_1\otimes\dots\otimes X_n)
\]
is $n$-linear, where each $\Tc X$ is a free $\Tc$-algebra. Moreover, it is universal in $\Cs^{\Tc,\varphi}$.
\end{lemma}
\begin{proof}
When $n=0$, the claim is about $\varphi_0:\unit\to\Tc\unit$. As mentioned, any map from $\unit$ to an algebra is $0$-linear. Moreover, $\varphi_0$ is universal by freeness of the algebra $\Tc\unit$.
When $n=1$, the claim is about $\id:\Tc X\to\Tc X$. But the identity of any algebra is $1$-linear and universal. For $n\geq 2$, we employ Lemma \ref{l:coeqmult}. We have
$\bar{\varphi}=\mu\circ\Tc\varphi$ and this is the coequalizer of the pair \eqref{eq:pair} by Lemma \ref{l:linton-free}. Hence $\varphi_{X_1,\dots,X_n}$ is multilinear and universal.
\end{proof}

We mention that the bilinearity of $\varphi_{X,Y}$ is in fact one of the conditions for a monoidal monad. The $n$-linearity of $\varphi_{X_1,\dots,X_n}$ may be deduced from this plus Lemma \ref{l:multisubcat}.

Together with Lemma \ref{l:post-tilde}, the following lemma addresses the behavior of the correspondence in Lemma \ref{l:coeqmult} with respect to composition. For $i=1,\dots,n$, let $f_i\colon X_i\to V_i$ be a morphism in $\Cs$, where $X_i$ is an object of $\Cs$ and $(V_i,v_i)$ is a $\Tc$-algebra.
Let $(W,w)$ be another algebra and $g: V_1\otimes \dots \otimes V_n \to W$ a $\Tc$-multilinear map.

\begin{lemma}\label{l:pre-tilde}
We have
\[
g \circ (\bar{f_1}\otimes \dots \otimes \bar{f_n})=\overline {g\circ (f_1\otimes \cdots \otimes f_n)} \circ \varphi.
\]
\end{lemma}
\begin{proof}
Since $g$ is multilinear and each $\bar{f_i}$ is a morphism of algebras, the composite on the left-hand side is multilinear. By universality of $\varphi$ (Lemma \ref{l:phi-multi}), there exists a (unique) map $h:X_1\otimes\dots\otimes X_n\to W$ such that
\[
g \circ (\bar{f_1}\otimes \dots \otimes \bar{f_n})= \bar{h}\circ \varphi.
\]
To determine $h$, precompose both sides above with $\eta_{X_1}\otimes\dots\otimes\eta_{X_n}$.
This yields
\[
g\circ (f_1\otimes\dots\otimes f_n) = \bar{h}\circ \eta_{X_1\otimes\dots\otimes X_n} = h,
\]
which completes the proof.
\end{proof}

Lemma \ref{l:post-tilde} is the case $n=1$ of Lemma \ref{l:pre-tilde}.

To a map $f:X_1\otimes\dots\otimes X_n\to Y$ in $\Cs$ we associate the composite
\[
\Tc X_1\otimes\dots\otimes\Tc X_n \map{\varphi} \Tc(X_1\otimes\dots\otimes X_n) \map{\Tc f} \Tc Y.
\]
Since $\varphi$ is $n$-linear (Lemma \ref{l:phi-multi}) and $\Tc f$ is $1$-linear, the composite is $n$-linear. This assignment is compatible with composition of multimaps. We obtain a morphism of multicategories
\[
\Mcc(\Cs,\otimes) \to \Cs^{\Tc,\varphi}
\]
which agrees with the free algebra functor $\Fc^\Tc:\Cs\to\Cs^\Tc$ on the underlying categories.

\subsection{Colax monoidal structure on the category of algebras}\label{ss:colax-mono-struct}

We are interested in potential monoidal structures on the category $\Cs^\Tc$, and in fact in colax monoidal structures  on this category. Recall from Section \ref{ss:multicat-colax} that any such structure $\und{\star}$ gives rise to a multicategory $\Mcc(\Cs^\Tc,\und{\star})$, and that if the structure is normal we have $\Uc\Mcc(\Cs^\Tc,\und{\star}) \cong\Cs^\Tc$.
The following proposition roughly states that the multicategories $\Mcc(\Cs^\Tc,\und{\star})$ and $\Cs^{\Tc,\varphi}$ are related provided that the (colax) monoidal categories $(\Cs^\Tc,\und{\star})$ and $(\Cs,\otimes)$ are related, and conversely.

To this end, we first regard the monoidal category $(\Cs,\otimes)$ as a normal colax monoidal category $(\Cs,\und{\otimes})$ (Section \ref{ss:colax}).

\begin{proposition}\label{p:multicat}
Fix a normal colax monoidal structure $\und{\star}$ on $\Cs^\Tc$.
There is a bijective correspondence between:
\begin{enumerate}[(i)]
\item \label{item:3} Monoidal structures $\gamma$ on the forgetul functor $\Uc^\Tc: (\Cs^\Tc,\und{\star})\to (\Cs, \und{\otimes})$ (between colax monoidal categories) for which the maps
\[
V_1\otimes\dots\otimes V_n = \otimes^n(\Uc^\Tc V_1,\dots,\Uc^\Tc V_n)\map{\gamma} \Uc^\Tc\bigl(\star^n(V_1 , \dots, V_n)\bigr) = \star^n(V_1 , \dots, V_n)
\]
are $\Tc$-multilinear;
\item \label{item:4}  Morphisms of multicategories
$\Fc\colon\Mcc(\Cs^\Tc,\und{\star})\to \Cs^{\Tc,\varphi}$
that restrict to the canonical isomorphism
\[
\Uc\Mcc(\Cs^\Tc,\und{\star})\cong \Cs^\Tc = \Uc(\Cs^{\Tc,\varphi})
\]
on the underlying categories.
\end{enumerate}
\end{proposition}
\begin{proof}
We sketch the argument.
Starting from $\gamma$ as in \eqref{item:3}, one defines $\Fc$ as follows. Let $f: (V_1,\dots,V_n)\to Y$ be a multimap in
$\Mcc(\Cs^\Tc,\und{\star})$. Thus, $f:\star^n(V_1,\dots,V_n)\to Y$ is a morphism in $\Cs^\Tc$.
Set
 \[
 \Fc(f) =f\circ \gamma.
 \]
Conversely, starting from $\Fc$ as in \eqref{item:4}), define $\gamma$ to be the image under $\Fc$
of the identity morphism of $\star^n(V_1,\dots,V_n)$ regarded as an $n$-map.
\end{proof}

We now construct a colax monoidal structure $\otimes^\varphi$ on $\Cs^\Tc$, under the assumption that
$\Cs$ possesses reflexive coequalizers and that these are preserved by $\Tc$.
The structure is such that $\Mcc(\Cs^\Tc,\otimes^\varphi)\cong \Cs^{\Tc,\varphi}$ as multicategories.

\begin{proposition}\label{p:universal}
Let $(\Tc, \varphi)$ be a monoidal monad on a monoidal category $(\Cs, \otimes)$.
Assume that $\Cs$ admits reflexive coequalizers and that these are preserved by $\Tc$. There
exists a normal colax monoidal structure $\und{\otimes}^\varphi$ on $\Cs^\Tc$ with the following properties.
 \begin{enumerate}[(i)]
 \item\label{item:7}
 The forgetful functor $\Uc^\Tc:(\Cs^\Tc,\und{\otimes}^\varphi) \to (\Cs,\und{\otimes})$ is monoidal
and its structure maps are $(\Tc, \varphi)$-multilinear.
 \item\label{item:8}
 The associated multicategory is the multicategory of  $(\Tc,\varphi)$-multilinear maps:
 \[
 \Mcc(\Cs^\Tc,\und{\otimes}^\varphi)\cong \Cs^{\Tc,\varphi}.
 \]
 \item\label{item:85}
 The preceding restricts to the canonical isomorphism $\Uc\Mcc(\Cs^\Tc,\und{\otimes}^\varphi)\cong \Cs^\Tc$ on the underlying categories.
 \end{enumerate}
\end{proposition}
Statements \eqref{item:8} and \eqref{item:85} say that the multicategory $\Cs^{\Tc,\varphi}$ is colax representable.
\begin{proof}
As discussed in Section \ref{ss:monad-Linton}, the hypotheses guarantee that Linton's coequalizers exist in $\Cs^\Tc$.
Lemma~\ref{l:coeqmult} then yields the
existence of universal $\Tc$-multilinear maps, and
Lemma~\ref{l:inducedmulticat} implies that $\Cs^{\Tc,\varphi}$ is colax
representable. Thus, there exists a normal colax monoidal structure $\und{\otimes}^\varphi$ on $\Cs^\Tc$ such that  $\Mcc(\Cs^\Tc,\und{\otimes}^\varphi)\cong \Cs^{\Tc,\varphi}$.
By Proposition~\ref{p:multicat}, the preceding isomorphism of multicategories corresponds to a monoidal structure on $\Uc:(\Cs^\Tc,\und{\otimes}^\varphi)\to(\Cs,\und{\otimes})$ whose components are $(\Tc,\varphi)$-multilinear.
%
\end{proof}

Proposition \ref{p:universal} holds under the weaker assumption on coequalizers in Lemma \ref{l:creation}, with the same proof.

We turn to the free algebra functor.

\begin{proposition}\label{p:universal-free}
Under the same hypothesis as for Proposition \ref{p:universal}, the free algebra functor
\[
\Fc^\Tc: (\Cs,\und{\otimes}) \to (\Cs^\Tc,\und{\otimes}^\varphi)
\]
is strong monoidal.
\end{proposition}
\begin{proof}
Universality of $\varphi_{X_1,\dots,X_n}:\Fc^\Tc X_1\otimes\dots\otimes\Fc^\Tc X_n\to \Fc^\Tc(X_1\otimes\dots\otimes X_n)$ in $\Cs^{\Tc,\varphi}$ (Lemma \ref{l:phi-multi}) affords canonical isomorphisms
\[
\otimes^{\varphi,n}(\Fc^\Tc X_1,\dots,\Fc^\Tc X_n) \cong \Fc^\Tc(X_1\otimes\dots\otimes X_n).\qedhere
\]
\end{proof}

The functor $\Fc^\Tc$ then induces a morphism of multicategories
\[
\Mcc(\Fc^\Tc): \Mcc(\Cs,\und{\otimes}) \to \Mcc(\Cs^\Tc,\und{\otimes}^\varphi) \cong \Cs^{\Tc,\varphi}.
\]
This agrees with the morphism in Section \ref{ss:free-multi}.

\subsection{Composition of universal multilinear maps}\label{ss:univ-multi}

We continue to work with a monoidal monad $(\Tc,\varphi)$ on $(\Cs, \otimes)$.

Consider algebras $(V_i,v_i)$ and $(V'_i,v'_i)$, $i=1,\dots,n$, $(W,w)$ and $(W',w')$. Consider a commutative diagram
\begin{equation}\label{eq:multilinear-comm}
\begin{gathered}
\xymatrix@C+20pt{
V_1\otimes  \dots \otimes V_n \ar[d]_h \ar[r]^-{f_1\otimes \dots \otimes f_n} & V'_1\otimes \dots \otimes V'_n \ar[d]^{h'}\\
W \ar[r]_-g & W'
}
\end{gathered}
\end{equation}
in which the $f_i$ and $g$ are in $\Cs^{\Tc}$, and $h$ and $h'$ are in $\Cs$.

\begin{lemma}\label{l:multilinear}
Suppose $h$ is $\Tc$-multilinear and $\Tc f_1\otimes \dots\otimes \Tc f_n$ is an epimorphism in $\Cs$.  Then $h'$ is $\Tc$-multilinear also.
\end{lemma}
\begin{proof}
Let $f=\Tc f_1\otimes \dots\otimes \Tc f_n$. We calculate:
\begin{align*}
w' \circ \Tc h' \circ \varphi \circ f & = w'\circ \Tc h' \circ \Tc(f_1\otimes \dots \otimes f_n) \circ \varphi & & \text{by naturality of $\varphi$,}\\
& = w' \circ \Tc g\circ \Tc h \circ \varphi & & \text{by commutativity,}\\
& = g \circ w \circ \Tc h \circ \varphi & & \text{since $g\in \Cs^\Tc$,}\\
& = g \circ h \circ (x_1\otimes \dots \otimes x_n) & & \text{by multilinearity of $h$,}\\
& = h' \circ (f_1\otimes \dots \otimes f_n) \circ  (x_1\otimes \dots \otimes x_n) & & \text{by commutativity,}\\
& = h' \circ (y_1\otimes \dots \otimes y_n) \circ f& & \text{since each $f_i \in \Cs^\Tc$.}
\end{align*}
Canceling $f$ we obtain the multilinearity of $h'$.
\end{proof}

Recall from~\cite[Lemma~4.7]{BW} and \cite[Lemma~0.17]{J}, that the tensor product functor
$\otimes: \Cs\times \Cs \to \Cs$ preserves reflexive coequalizers (as a functor of two variables) if and only if for each object $X$ in $\Cs$, the one-variable functors $X\otimes(-)$ and $(-)\otimes X$ preserve reflexive coequalizers. We make this assumption below.

Let $X$ be an object of $\Cs$ and $A$ a $\Tc$-algebra. We say that a map $e\colon X\to A$
in $\Cs$ is a \emph{reflexive $\Tc$-coequalizer} if the corresponding morphism
$\bar e\colon \Tc X\to A$ is a reflexive coequalizer in $\Cs^\Tc$. 

Universal multilinear maps are examples of reflexive $\Tc$-coequalizers, by \eqref{item:6} in Lemma \ref{l:coeqmult}. Such maps need not be epimorphisms in $\Cs$: the universal bilinear map for a pair of free algebras fails to be surjective for the multiset monad and for the powerset monad (Sections \ref{ss:multiset} and \ref{ss:powerset}).
However, reflexive $\Tc$-coequalizers possess the following property.

\begin{lemma}\label{l:univ-epi}
Let  $u_i: X_i\to A_i$ be a reflexive $\Tc$-coequalizer, $i=1,\dots,n$. Let $f$ and $g: A_1\otimes\dots\otimes A_n\to B$ be two $\Tc$-multilinear maps. Assume that $\Cs$ admits reflexive coequalizers and both $\otimes$ and $\Tc$ preserve such coequalizers. If 
\[
f\circ (u_1\otimes\dots\otimes u_n)=g\circ (u_1\otimes\dots\otimes u_n),
\] 
then in fact $f=g$.
\end{lemma}
\begin{proof}
By Lemma \ref{l:pre-tilde},
\[
f\circ (\overline{u_1}\otimes\dots\otimes \overline{u_n})= 
\overline{f\circ(u_1\otimes\dots\otimes u_n)}\circ\varphi =
\overline{g\circ(u_1\otimes\dots\otimes u_n)}\circ\varphi =
g\circ (\overline{u_1}\otimes\dots\otimes \overline{u_n}).
\]
The map $\overline{u_1}\otimes\dots\otimes \overline{u_n}$ is a reflexive coequalizer in $\Cs$ by the assumptions. It follows that $f=g$.
\end{proof}

\begin{proposition}\label{p:main}
Assume that $\Cs$ admits reflexive coequalizers and both $\otimes$ and $\Tc$ preserve such coequalizers. Then, universal $\Tc$-multilinear maps are closed under composition.
\end{proposition}
\begin{proof}
For each $i=1,\dots,n$, let $\und{U^i}$ be a sequence of algebras and let $U_i$ denote the tensor product of the terms in each sequence. Let
$a_i:U_i\to V_i$ be a universal multilinear map for $\und{U^i}$.
Let $b:V_1\otimes\dots\otimes V_n\to W$ be an additional universal multilinear map.
In order to show that the composite
\[
U_1\otimes\dots\otimes U_n \map{a_1\otimes\dots\otimes a_n} V_1\otimes\dots\otimes V_n \map{b} W
\]
is universal, we need to verify that any multilinear map 
\[
g:U_1\otimes\dots\otimes U_n \to X
\]
factors uniquely as the composite above followed by a morphism of algebras. We construct the factorization in two steps.

Consider the sequence obtained by applying $\Tc$ to each term in the sequence
$\und{U^i}$. Let $T_i$ denote the tensor product of the terms in this sequence. Let $u_i:T_i\to U_i$ denote the tensor product of the structure maps of the algebras in the sequence $\und{U^i}$.
We claim that we also have two commutative squares as follows.
\[
\xymatrix@C+50pt{
\Tc T_1\otimes  \dots \otimes \Tc T_n \ar[d]_{\varphi}  \ar@<0.5ex>[r]^-{(\mu\circ\Tc\varphi)\otimes\dots\otimes(\mu\circ\Tc\varphi)} \ar@<-0.5ex>[r]_-{\Tc u_1\otimes\dots\otimes \Tc u_n}  & \Tc U_1\otimes \dots \otimes \Tc U_n \ar[d]^{\varphi}\\
\Tc(T_1\otimes  \dots \otimes T_n) \ar@<0.5ex>[r]^-{\mu\circ\Tc\varphi} \ar@<-0.5ex>[r]_-{\Tc(u_1\otimes\dots u_n)} & \Tc(U_1\otimes \dots \otimes U_n)
}
\]
The top row is the tensor product of the Linton pairs associated to the sequences $\und{U^i}$.
The bottom row is the Linton pair associated to their concatenation. The bottom squares commutes by naturality of $\varphi$. The commutativity of the top square follows from associativity for the monoidal structure $\varphi$ and multilinearity of $\varphi$ (Lemma \ref{l:phi-multi}).

By item \eqref{item:5} in Lemma \ref{l:coeqmult}, $\bar{g}$ coequalizes the bottom pair. 
By item \eqref{item:6} in the same result, and since $\otimes$ preserves coequalizers, 
$\overline{a_1}\otimes\dots\otimes\overline{a_n}$ is the coequalizer of the top pair.
We deduce a (unique) map $\tilde{g}$ in $\Cs$ as below
\[
\xymatrix@C+20pt{
\Tc U_1\otimes\dots\otimes\Tc U_n \ar[r]^-{\overline{a_1}\otimes\dots\otimes\overline{a_n}} \ar[d]_{\varphi} &
 V_1\otimes\dots\otimes V_n \ar@{-->}[d]^{\tilde{g}} \\
 \Tc(U_1\otimes\dots\otimes U_n) \ar[r]_-{\bar{g}} & X.
}
\]
By Lemma \ref{l:multilinear}, $\tilde{g}$ is multilinear. (Note that $\varphi$ is multilinear by Lemma \ref{l:phi-multi} and $\Tc\overline{a_1}\otimes\dots\otimes\Tc\overline{a_n}$ is an epimorphism since $\Tc$ and $\otimes$ preserve coequalizers.) Then, universality of $b$ yields a (unique) morphism of algebras $\hat{g}$ as below.
\[
\xymatrix{
 V_1\otimes\dots\otimes V_n \ar[d]_{\tilde{g}} \ar[r]^-{b} & W \ar@{-->}[dl]^{\hat{g}}\\
 X & 
}
\]
By construction,
\[
\hat{g}\circ b\circ(\overline{a_1}\otimes\dots\otimes\overline{a_n}) = \bar{g}\circ\varphi.
\]
Precomposing with $\eta_{V_1}\otimes\dots\otimes\eta_{V_n}$ we obtain
\[
\hat{g}\circ b\circ(a_1\otimes\dots\otimes a_n) = g.
\]
This is the desired factorization. Its uniqueness follows from the uniqueness in each of the two steps of the construction, or as an application of Lemma \ref{l:univ-epi}.
\end{proof}

\subsection{Monoidal structure on the category of algebras}\label{ss:monoidal}

We arrive at the main result of the section: the existence of a monoidal structure on the category of algebras over a monoidal monad. We show that the colax monoidal structure on $\Cs^\Tc$ from Proposition \ref{p:universal} is in fact monoidal. This builds on the work in Section \ref{ss:univ-multi} and the hypothesis are as in Proposition \ref{p:main}.

\begin{theorem}\label{t:monoidal}
Let $(\Tc,\varphi)$ be a monoidal monad on the monoidal category $(\Cs, \otimes)$.
Assume that $\Cs$ admits reflexive coequalizers and that both $\otimes$ and $\Tc$ preserve such coequalizers. Then the colax monoidal structure $\und{\otimes}^\varphi$ on $\Cs^\Tc$ is in fact monoidal.
\end{theorem}
\begin{proof}
According to Lemma \ref{l:inducedmulticat}, the colax monoidal structure is monoidal when universal multilinear maps are closed under composition. This is the result of Proposition \ref{p:main}.
\end{proof}

We denote the tensor product on $\Cs^\Tc$ by $\otimes^\varphi$.
The following recapitulates some of the results from previous sections.
The tensor product of two algebras $(V,v)$ and $(W,w)$ is the Linton coequalizer
\[
\xymatrix{
	\Tc(\Tc V \otimes \Tc W) \ar[r]^-{\Tc\varphi}
	\ar@/_1pc/[rr]_-{\Tc(v\otimes w)} &
	\Tc^2(V \otimes W) \ar[r]^-{\mu}
	& \Tc(V \otimes W) \ar[r] & V\otimes^\varphi W
}
\]
in $\Cs^\Tc$.
The unit object is $(\Tc\unit, \mu_\unit)$.
The free algebra functor 
\[
\Fc^\Tc: (\Cs,\otimes,\unit) \to (\Cs^\Tc,\otimes^\varphi,\Tc\unit)
\] 
is strong monoidal.

\subsection{Example: monads on cocartesian categories}\label{ss:cocartesian}

Let $\Cs$ be a category admitting all finite coproducts. We may view it as a monoidal category in which the tensor product is the categorical coproduct $\amalg$ and the unit object is an initial object $\bot$. Categories of this form are called \emph{cocartesian}.

Any functor between cocartesian categories carries a unique monoidal structure. (See, for instance, \cite[Example 3.19]{AM}.) In addition, this structure is preserved by any natural transformation between two such functors. It follows that any monad $\Tc$ on $\Cs$ admits a unique monoidal structure. Explicitly, the structure is
\[
\varphi_{X,Y}: \Tc X\amalg \Tc Y \map{\binom{\Tc\iota_X}{\Tc\iota_Y}} \Tc(X\amalg Y),
\]
where $\iota_X:X\to X\amalg Y$ and $\iota_Y:Y\to X\amalg Y$ are the canonical insertions.

Linton studied the existence of coproducts in the category of algebras over such a monad.
He showed that coproducts in $\Cs^\Tc$ are defined by the coequalizers \eqref{eq:pair} \cite[Proposition 2]{Li}.
We review his result from our perspective.

Let $U$, $V$ and $W$ be $\Tc$-algebras. One easily verifies that a map $\binom{f_1}{f_2}:U\amalg V\to W$ is $\Tc$-bilinear if and only if  its components $f_1:U\to W$ and $f_2:V\to W$ are morphisms of $\Tc$-algebras. It follows that the coproduct of $U$ and $V$ exists in $\Cs^\Tc$ precisely when the universal $\Tc$-bilinear map exists for $(U,V)$. In this case, the coproduct $U\amalg^\varphi V$ in $\Cs^\Tc$ is given by the Linton coequalizer
\[
\xymatrix{
	\Tc(\Tc U \amalg \Tc V) \ar[r]^-{\Tc\varphi}
	\ar@/_1pc/[rr]_-{\Tc(u\amalg v)} &
	\Tc^2(U \amalg V) \ar[r]^-{\mu}
	& \Tc(U \amalg V) \ar[r] & U\amalg^\varphi V
}
\]
in $\Cs^\Tc$.
It also follows that the coinsertions in the coproduct $U\amalg^\varphi V$ are the components of the universal bilinear map $U\amalg V\to U\amalg^\varphi V$.

Since the resulting structure on $\Cs^\Tc$ is the cocartesian one,
the fact that this structure is monoidal is guaranteed with no assumptions beyond existence of  Linton coequalizers.

\subsection{Example: linear monads}\label{ss:linear-monoidal-monad}

Let $(\Cs,\otimes)$ be a braided monoidal category with braiding $\beta$.
Let $A$ be a commutative monoid in $\Cs$. The linear monad $\Tc=A\otimes(-)$ is monoidal on $(\Cs,\otimes)$. The structure $\varphi_{X,Y}$ is
\[
A\otimes X\otimes A\otimes Y \map{\id\otimes\beta\otimes\id} A\otimes A\otimes X\otimes Y \map{\mu\otimes\id}A\otimes X\otimes Y.
\]
Commutativity of $A$ guarantees the compatibility between the monoidal and the monad structures.

Let $L$, $M$ and $N$ be left $A$-modules in $\Cs$. Let $m$ and $n$ denote the actions of $A$ on $M$ and $N$. Kock's criterion (Lemma \ref{l:kock}) says that a map $f:M\otimes N\to L$ is $\Tc$-bilinear precisely when it satisfies the following two conditions.
\begin{itemize}
\item $f$ is a morphism of $A$-modules with $M\otimes N$ viewed as an $A$-module under
\[
A\otimes M\otimes N \map{m\otimes\id} M\otimes N;
\]
\item $f$ is a morphism of $A$-modules with $M\otimes N$ viewed as an $A$-module under
\[
A\otimes M\otimes N \map{\beta\otimes\id}	
M\otimes A\otimes N \map{\id\otimes n}  M\otimes N.
\]
\end{itemize}
It follows that
\[
M\otimes^\varphi N = M\otimes_A N,
\]
the familiar tensor product of modules over a commutative monoid.


The preceding example may be generalized: starting from a \emph{double monoid} $A$ in a duoidal category $(\Cs,\diamond,\star)$, the monad $H\star(-)$ is monoidal on the duoidal category $(\Cs,\diamond)$. Note the order of the operations. See \cite[Section 6.5]{AM} for information on double monoids.

\subsection{Example: the multiset monad}\label{ss:multiset}

Let $\Nb$ denote the set of natural numbers. It is a commutative semiring. 
Let $\Set$ denote the category of sets. 

Given a set $X$, let $\Tc X$ be the set of functions
\[
f:X\to\Nb
\]
of finite support. Such a function is a \emph{multisubset} of $X$. Given a map $h:X\to Y$, let $\Tc h:\Tc X\to\Tc Y$ be defined by
\[
(\Tc h)(f)(y) = \sum_{x\in h^{-1}(y)} f(x).
\]
The endofunctor $\Tc$ on $\Set$ carries a monad structure. The multiplication $\mu$ is as follows. Let $F$ be an element of $\Tc^2 X$. It is a function $\Tc X\to\Nb$ of finite support. Then $\mu_X(F)\in\Tc X$ is the function
\[
X\to\Nb, \quad x \mapsto \sum_{f\in\Tc X} F(f) f(x).
\]
It is well-defined since $F$ is of finite support. The unit $\eta$ is as follows: $\eta_X(x)\in\Tc X$ is the function
\[
X\to\Nb, \quad y\mapsto \delta(x,y)
\]
(the Kronecker delta). Algebras over $(\Tc,\mu,\eta)$ are commutative monoids in the standard sense (commutative monoids in the cartesian category $(\Set,\times,\top)$): such a monoid $A$ is a $\Tc$-algebra via 
\[
\alpha:\Tc A\to A, \quad \alpha(f)=\sum_{a\in A} f(a)\cdot a.
\]

The monad $\Tc$ on $(\Set,\times,\top)$ carries a monoidal structure. The structure transformation $\varphi$ is as follows. Let $(f,g)$ be an element of $\Tc X\times\Tc Y$. Then $\varphi_{X,Y}(f,g)\in\Tc(X\times Y)$ is the function
\[
X\times Y\to\Nb, \quad (x,y)\mapsto f(x)g(y).
\]
The map $\varphi_0$ sends the unique element of $\top$ to $1\in\Nb=\Tc\top$.

The monad $\Tc$ is associated to the \emph{algebraic theory} of commutative monoids.
It follows that is preserves reflexive coequalizers \cite[Corollaries A.22 and A.23]{ARV2011}.
(This also follows by a direct argument.) 
Theorem \ref{t:monoidal} yields a monoidal structure $\otimes^\varphi$ on the category of commutative monoids.

Let $A$, $B$ and $C$ be commutative monoids. Kock's criterion (Lemma \ref{l:kock}) says that a map $h:A\times B\to C$ is $\Tc$-bilinear if and only if $h(-,b):A\to C$ and $h(a,-):A\to C$ are morphisms of monoids for each $a\in A$ and $b\in B$.
It follows that $A\otimes^\varphi B$ is the familiar tensor product of commutative monoids.

Abelian groups and modules over commutative rings may be treated similarly.

\subsection{Example: the powerset monad} \label{ss:powerset}

The preceding example may be generalized: starting from a \emph{commutative semiring} $R$ (with additive and multiplicative units), there is a monoidal monad $\Tc$ on $(\Set,\times,\top)$ defined by 
\[
\Tc X = \{f:X\to R \mid \text{$f$ is of finite support}\}.
\]
This leads to a notion of tensor product for \emph{$R$-semimodules} \cite[Chapter 16]{Gol1999}.

Consider the case in which $R=\{0,1\}$ is the Boolean semiring. In this case $\Tc X$ may be identified with the collection of finite subsets of $X$, and $\Tc$ is the finite \emph{powerset monad}.
Algebras are \emph{$\vee$-semilattices}: sets equipped with a binary operation $\vee$ that is associative, unital, commutative, and idempotent.
For the same reasons as above, this monad preserves reflexive coequalizers. We obtain a monoidal structure on the category of $\vee$-semilattices.

Consider instead the full powerset monad, for which $\Tc X$ is the collection of all subsets of $X$. Algebras are now \emph{complete $\vee$-semilattices}. This monad does not correspond to an algebraic theory and in fact does not preserve reflexive coequalizers. Theorem \ref{t:monoidal} does not apply. Nevertheless, the resulting structure on the category of complete $\vee$-semilattices is monoidal.
This is guaranteed by an alternative set of hypotheses for Theorem \ref{t:monoidal} which we do not discuss in this paper (the closedness of the ground category).

The tensor product of semilattices is considered in \cite{Fra1976} and \cite[Section I.5]{JT}.

Here is a reflexive pair whose coequalizer is not preserved by the full powerset monad. Let $f,g:\Nb \times \{0,1\} \to\Nb$ be given by 
\[
f(n,0)=f(n,1)=n \qand g(n,0)=n,\ g(n,1)=n+1.
\] 
The function $s(n)=(n,0)$ is a common section. See \cite[Proposition 4.6.5]{Bor:1994ii} for complementary information on the powerset monad.

\section{The pointing monad and the smash product}\label{s:pointing}

We illustrate the constructions of the preceding sections with the pointing monad and related notions. The goal is to arrive at the definition of smash product of bipointed objects in a monoidal category.

\subsection{Pointed, copointed, and bipointed objects}\label{ss:pointed}

Let $(\Cs,\unit)$ be a \emph{pointed} category: a category together with a distinguished object $\unit$.  A \emph{pointed object} is a pair $(X,e)$ where $e:\unit\to X$ is a map in $\Cs$.  The map $e$ is the \emph{point} of $X$. The morphisms of pointed objects are the point-preserving maps in $\Cs$. Let $\Csp$ denote the category of pointed objects in $(\Cs,\unit)$.
It is the coslice category of objects under $\unit$.

A \emph{copointed object} is a pair $(X,\epsilon)$ where $\epsilon:X\to\unit$ is a map in $\Cs$.
Let $\Csc$ denote the category of copointed objects. It is the slice category of objects over $\unit$.

An object $X$ is \emph{bipointed} if it is equipped with a point $e:\unit\to X$ and a copoint $\epsilon:X\to\unita$ such that
\[
\xymatrix@-10pt{
 & X \ar[rd]^{\epsilon} & \\
\unit \ar[ru]^{e} \ar@{=}[rr]_{} & & \unit
}
\]
commutes. Let $\Cspc$ denote the category of bipointed objects in $\Cs$.

The object $\unit$ is bipointed with the point and copoint equal to $\id_\unit$. As a pointed object, it is initial in the category $\Csp$. As a copointed object, it is terminal in $\Csc$. As a bipointed object, it is null in $\Cspc$.

Suppose the object $\unit$ is terminal in $\Cs$. In this case, every object of $\Cs$ carries a unique copoint, and there are canonical equivalences
\[
\Cs \cong \Csc \qand \Csp \cong \Cspc.
\]
There are corresponding statements when $\unit$ is initial or null.

\subsection{The smash product of bipointed objects}\label{ss:smash-pushout}

Let $(\Cs,\otimes,\unit)$ be a monoidal category. Consider the pointed category $(\Cs,\otimes)$. Assume $\Cs$ admits finite coproducts.
The \emph{smash product} $X\wedge Y$ of two bipointed objects $X$ and $Y$ is defined as the following pushout in the category $\Csc$ (when this exists):
\begin{equation}\label{eq:smash}
\begin{gathered}
\xymatrix@C+10pt{
X\amalg Y \ar[r]^-{\binom{\id_X\otimes e_Y}{e_X\otimes\id_Y}} \ar[d]_{\binom{\epsilon_X}{\epsilon_Y}} & X\otimes Y \ar[d] \\
\unit \ar[r] & X\wedge Y.
}
\end{gathered}
\end{equation}
The components of the top map are
\[
\xymatrix@C+10pt{
X=X\otimes\unit \ar[r]^-{\id_X\otimes e_Y} & X\otimes Y & \unit\otimes Y=Y. \ar[l]_-{e_X\otimes\id_Y}
}
\]
The bottom map defines the point of $X\wedge Y$.

Anel and Joyal have considered this notion in \cite[Section 1.1.1]{AJ}.

In the following sections we discuss how these notions can be formulated in terms of monads and the constructions of the previous sections.



\subsection{The pointing monad}\label{ss:pointing-monad}

We return to a pointed category $(\Cs,\unit)$.
Assume $\Cs$ admits finite coproducts. The \emph{pointing monad} $\Tcp$ is defined by
\[
\Tcp X = \unit\amalg X.
\]
Algebras over $\Tcp$ are precisely pointed objects in $(\Cs,\unit)$. It is a linear monad as in Section \ref{ss:linear-monad}, since $\unit$ is a monoid in $(\Cs,\amalg,\bot)$.

Assume that $(\Cs,\otimes,\unita)$ is monoidal. Suppose that the object $\unit$ is a comonoid in $\Cs$, but not necessarily the unit object $\unita$. In this case, the monad $\Tcp$ is comonoidal. Indeed, we may view $\unit$ as a bimonoid in the duoidal category $(\Cs,\amalg,\bot,\otimes,\unita)$ \cite[Example 6.19]{AM}, and the construction of Section \ref{ss:linear-monad} applies. According to Theorem \ref{t:comonoidal}, the monoidal structure lifts to the category of $\Tcp$-coalgebras. Explicitly, if $X$ and $Y$ are pointed, then so is $X\otimes Y$ with point
\[
\unit\map{}\unit\otimes\unit \map{e_X\otimes e_Y} X\otimes Y.
\]
The unlabeled map is the coproduct of $\unit$. The point of $\unita$ is the counit of $\unit$. We obtain a monoidal category $(\Csp,\otimes,\unita)$.

We turn to a monoidal structure on $\Tcp$. It is defined under a different set of assumptions from the preceding.
Suppose the functors $X\otimes(-)$ and $(-)\otimes X$ preserve coproducts, for each object $X$ of $\Cs$. Again $\unit$ need not be the unit object of $\Cs$. Instead, suppose $\unit$ is terminal in $\Cs$. We may then define a map as follows, for any objects $X$ and $Y$ of $\Cs$:
\[
(\unit\amalg X)\otimes(\unit\amalg Y) \map{\cong} (\unit\otimes\unit)\amalg(\unit\otimes Y)\amalg(X\otimes\unit)\amalg(X\otimes Y) \to \unit\amalg(X\otimes Y).
\]
Preservation of coproducts affords the isomorphism on the left. The map on the right has the following components:
\[
\unit\otimes\unit \to \unit, \quad
\unit\otimes Y \to  \unit, \quad
X\otimes\unit \to \unit, \quad
X\otimes Y \map{\id} X\otimes Y.
\]
The unlabeled maps are unique.

\begin{lemma}\label{l:pointing-monoidal}
With the above structure, the monad $\Tcp$ is monoidal on $(\Cs,\otimes)$.
\end{lemma}

Note that even though the monad $\Tcp$ is linear, the underlying monoidal structure is not, in the sense that it does not arise from the construction of Section \ref{ss:linear-monoidal-monad}. For this, one would require first of all a duoidal structure for $(\Cs,\otimes,\amalg)$ (and then a double monoid structure on $\unit$). But in general there is no such structure.

\subsection{The copointing comonad}\label{ss:copointing-comonad}

Assume now that $\Cs$ admits finite products.
We denote them by $X\times Y$ and use $\pi$ for the canonical projections.
The \emph{copointing comonad} $\Tcc$ is defined by
\[
\Tcc X = \unit\times X.
\]
Coalgebras over $\Tcc$ are precisely copointed objects in $\Cs$.

Assume $\Cs$ admits finite coproducts as well as finite products.
There is a \emph{mixed distributive law} $\Tcp\Tcc \to \Tcc\Tcp$ defined as follows. The map
\[
\unit\amalg(\unit\times X) \to \unit\times(\unit\amalg X)
\]
has components
\[
\per{\id_\unit & \iota_\unit \\ \pi_\unit & \iota_X\circ \pi_X}.
\]
The bottom right corner is the composite $\unit\times X\map{\pi_X} X \map{\iota_X}\unit\amalg X$ between a product projection and a coproduct insertion. It follows that the monad $\Tcp$ lifts to the category of copointed objects. Explicitly, if $\epsilon$ is the copoint of $X$ , then
\[
\unit\amalg X \map{\binom{u}{\epsilon}} \unit
\]
is the copoint of $\Tcp X$. Algebras over the lifted monad are precisely bipointed objects.
We continue to write $\Tcp$ for the monad lifted to $\Csc$.

\subsection{The smash product as a monoidal structure}\label{ss:smash}

Let $(\Cs,\otimes,\unit)$ be a monoidal category. We view it as a pointed category $(\Cs,\unit)$ and consider (co)pointed objects in this category, and the associated monad and comonad.

Assume that $\Cs$ admits finite products, so the comonad $\Tcc$ is defined.
The comonad $\Tcc$ is monoidal on $\Cs$. In fact, the unit object $\unit$ is a monoid in this category, so the situation is dual to that in Section \ref{ss:pointing-monad} (where we discussed the comonoidal structure of $\Tcp$). Hence, the monoidal structure lifts to the category of $\Tcc$-coalgebras. Explicitly, if $X$ and $Y$ are copointed, then so is $X\otimes Y$ with copoint
\[
X\otimes Y \map{\epsilon_X\otimes\epsilon_Y} \unit\otimes \unit \cong\unit.
\]
The copoint of $\unit$ is $\id_\unit$. We obtain a monoidal category $(\Csc,\otimes,\unit)$.

Assume that $\Cs$ admits finite coproducts and that they are preserved
by the functors $X\otimes(-)$ and $(-)\otimes X$, for each object $X$.
Since the forgetful functor $\Csc\to\Cs$ creates colimits, the monoidal category $\Csc$ inherits these properties. It follows that the lifted monad $\Tcp$ on $\Csc$ is the pointing monad associated to the terminal object $\unit$ of $\Csc$.
Lemma \ref{l:pointing-monoidal} then affords a monoidal structure for the monad $\Tcp$ on $(\Csc,\otimes,\unit)$.


Let $X$, $Y$ and $Z$ be bipointed objects. A morphism $f:X\otimes Y\to Z$ of copointed objects is $\Tcp$-bilinear precisely when the diagram
\begin{equation}\label{eq:smash-bilinear}
\begin{gathered}
\xymatrix@C+10pt{
X\amalg Y \ar[r]^-{\binom{\id_X\otimes e_Y}{e_X\otimes\id_Y}} \ar[d]_{\binom{\epsilon_X}{\epsilon_Y}} & X\otimes Y \ar[d]^{f} \\
\unit \ar[r]_{e_Z} & Z
}
\end{gathered}
\end{equation}
commutes. It follows that the universal $\Tcp$-bilinear map for $(X,Y)$ exists if and only if the smash product $X\wedge Y$ exists. In this case, the right vertical arrow in \eqref{eq:smash} is the universal map of $(X,Y)$, and the bottom horizontal arrow is the point of $X\wedge Y$.

Thus, the smash product $X\wedge Y$ is the colax monoidal structure afforded by Proposition \ref{p:universal}. The situation for $n$ copointed objects is similar.

Theorem \ref{t:monoidal} provides conditions that guarantee existence and associativity of the smash product. Note that the lifted monad $\Tcp$ preserves colimits. Therefore, to obtain the preceding, it suffices to assume that $\Csc$ admits finite colimits and that these are preserved by the one-variable tensor functors.

\subsection{The smash product of pointed comonoids}\label{ss:smash-comonoid}

Let $\Co{\Cs}$ denote the category of comonoids in $(\Cs,\otimes,\unit)$. An object is a triple $(C,\Delta,\epsilon)$ subject to the usual axioms. Forgetting the coproduct $\Delta$ and retaining the counit $\epsilon$ as the copoint yields a functor
\[
\Co{\Cs} \to \Csc.
\]

We may regard $\unit$ as a bimonoid in the duoidal category $(\Cs,\amalg,\bot,\otimes,\unit)$, and therefore $\Tcp$ as a comonoidal monad on $(\Cs,\otimes,\unit)$. It follows that $\Tcp$ lifts to $\Co{\Cs}$ \cite[Proposition 2.1]{M}. Algebras over the lifted monad are pointed comonoids in $\Cs$.

Assume now that $(\Cs,\otimes,\unit)$ is braided. The category $\Co{\Cs}$ is then monoidal, and the monad $\Tcp$ on this category is monoidal in the same manner as in Lemma \ref{l:pointing-monoidal}. This leads to the notion of smash product of pointed comonoids. It is defined by the same diagram \eqref{eq:smash}, where the pushout is now calculated in the category of comonoids. This notion is considered in \cite[Section 1.3.6]{AJ}.

\section{Bimonoidal monads}\label{s:bimonoidal}

Consider a monad $\Tc$ acting on a duoidal category $(\Cs, \diamond, \star)$. The monad may interact with the monoidal structures in three different ways. Here we treat the case in which $\Tc$ is monoidal with respect to $\diamond$ and comonoidal with respect to $\star$; the other two cases are treated in the following sections.
Let $\varphi$ and $\psi$ denote the monoidal and comonoidal structures of $\Tc$, respectively.
The (colax) monoidal structure on $\Cs^\Tc$ resulting from the former is denoted by $\diamond^\varphi$. It is preserved by the free algebra functor. This was the subject of Section \ref{s:monoidal}. The monoidal structure resulting from the latter is still denoted by $\star$. It is preserved by the forgetful functor. This was the subject of Section \ref{s:comonoidal}. Our goal is to show that $(\Cs^\Tc,\diamond^\varphi,\star)$ is itself a duoidal category. This requires a compatibility between the monoidal and comonoidal structures on the monad.

\subsection{Bimonoidal monads}

We say that $(\Tc, \varphi, \psi)$ is a \emph{bimonoidal} monad if (in addition to the above conditions) the following diagrams commute.
\begin{equation}\label{eq:bilaxass}
\begin{gathered}
\xymatrix@C=-40pt@R=+20pt{
& \Tc(A \star B) \diamond \Tc(C \star D)
\ar[rd]^-{\psi \diamond \psi}
\ar[ld]_-{\varphi} & \\
{\Tc\bigl((A \star B) \diamond (C \star D)\bigr)} \ar[d]_{\Tc \zeta}  & &
(\Tc A \star \Tc B) \diamond (\Tc C \star \Tc D)
\ar[d]^{\zeta}\\
{\Tc\bigl((A \diamond C) \star (B \diamond D)\bigr)} \ar[rd]_-{\psi} & &
(\Tc A \diamond \Tc C) \star (\Tc B \diamond \Tc D)
\ar[dl]^-{\varphi \star \varphi} \\
& \Tc(A \diamond C) \star \Tc(B \diamond D)
}
\end{gathered}
\end{equation}

\begin{equation}\label{eq:ibunit}
\begin{gathered}
\xymatrix@R-0.5pc@C-.5pc{
\unit \ar[d]_{}\ar[r]^-{\varphi_0} &
\Tc \unit \ar[r]^-{} &
\Tc(\unit \star \unit) \ar[d]^{\psi_{\unit,\unit}}\\
\unit \star \unit \ar[rr]_-{\varphi_0 \star \varphi_0} & &
\Tc \unit  \star \Tc \unit
}
\qquad
\xymatrix@R-0.5pc@C-.5pc{
\unita & \Tc\unita \ar[l]_-{\psi_0} &
\Tc(\unita \diamond \unita) \ar[l]_-{}\\
\unita \diamond \unita \ar[u]^{} & &
\Tc \unita \diamond \Tc\unita
\ar[ll]^-{\psi_0 \diamond \psi_0} \ar[u]_{\varphi_{\unita,\unita}}
}
\end{gathered}
\end{equation}

\begin{equation}\label{eq:ibinverse}
\begin{gathered}
\xymatrix@R-0.5pc@C-1pc{
\Tc \unit \ar[r]^{} &
\Tc \unita \ar[d]^-{\psi_0}\\
\unit \ar[r]_{}\ar[u]^-{\varphi_0} & \unita
}
\end{gathered}
\end{equation}
The unlabeled arrows are either unit maps of $\Cs$ or their images under $\Tc$.

Consider the $2$-category whose objects are duoidal categories, whose morphisms are bilax monoidal functors, and whose $2$-cells are morphisms of such functors \cite[Proposition 6.52]{AM}. A bimonoidal monad is precisely a monad in this $2$-category.

We remark that \eqref{eq:ibinverse} and the first diagram in \eqref{eq:ibunit} are superfluous: the commutativity of these diagrams follows from the fact that $\eta$ is a morphism of monoidal functors and of comonoidal functors, which includes the properties
\[
\varphi_0 = \eta_\unit \qand \psi_0\circ\eta_\unita = \id_{\unita}.
\]
The two remaining conditions may be interpreted in terms of $(\Tc,\varphi)$-bilinearity for the monoidal structures on algebras induced by $\psi$, as we next discuss.

The comonoidal monad $(\Tc,\psi)$ on $(\Cs,\star,\unita)$ gives rise to a monoidal structure  on $\Cs^\Tc$ preserved by the forgetful functor (Theorem \ref{t:comonoidal}). Given $\Tc$-algebras $X$ and $Y$, $X\star Y$ is a $\Tc$-algebra via
\[
\Tc(X\star Y)\map{\psi_{X,Y}} \Tc X\star\Tc Y\map{x\star y} X\star Y.
\]
The unit object $\unita$ is a $\Tc$-algebra via $\Tc\unita\map{\psi} \unita$. We do not distinguish notationally between the tensor product $\star$ on $\Cs$ and on $\Cs^\Tc$.

The second diagram in \eqref{eq:ibunit} says that the unit map $\unita\diamond\unita\to\unita$ is $(\Tc,\varphi)$-bilinear. Regarding \eqref{eq:bilaxass}, we have the following result.

\begin{lemma}\label{l:bimonoidal-duoidal}
Let $A$, $B$, $C$, $D$, $X$ and $Y$ be $\Tc$-algebras. Let $f:A\diamond C\to X$ and $g:B\diamond D\to Z$ be $(\Tc,\varphi)$-bilinear maps. Then
\[
(A\star B)\diamond(C\star D) \map{\zeta} (A\diamond C)\star(B\diamond D) \map{f\star g} X\star Y
\]
is $(\Tc,\varphi)$-bilinear too.
\end{lemma}
\begin{proof}
Bilinearity is shown by the commutativity of the diagram below.
\[
\def\objectstyle{\scriptstyle}
\def\labelstyle{\scriptstyle}
\xymatrix@C=-5pt{
& \Tc(A \star B) \diamond \Tc(C \star D)
\ar[rd]^-{\psi \diamond \psi}
\ar[ld]_-{\varphi} & \\
{\Tc(A \star B) \diamond (C \star D))} \ar[d]_{\Tc \zeta}  & &
(\Tc A \star \Tc B) \diamond (\Tc C \star \Tc D) \ar[d]^{\zeta} \ar[rd]^-{(a\ast b)\diamond(c\ast d)} \\
{\Tc((A \diamond C) \star (B \diamond D))} \ar[rd]_-{\psi} \ar[d]_-{\Tc(f\ast g)} & &
(\Tc A \diamond \Tc C) \star (\Tc B \diamond \Tc D)
\ar[dl]^-{\varphi \star \varphi} \ar[rd]_-{(a\diamond c)\star(b\diamond d)} & (A \star B) \diamond (C \star D) \ar[d]^-{\zeta} \\
\Tc(X\ast Y) \ar[dr]_-{\psi} & \Tc(A \diamond C) \star \Tc(B \diamond D) \ar[d]^{\Tc f\star\Tc g} & & (A \diamond C) \star (B \diamond D) \ar[d]^-{f\ast g}\\
& \Tc X\ast\Tc Y \ar[rr]_-{x\ast y} & & X\ast Y
}
\]
The hexagon is axiom \eqref{eq:bilaxass} for the bimonoidal monad and the pentagon commutes by bilinearity \eqref{eq:multilinear} of $f$ and $g$.
\end{proof}

Bimonoidal monads are in fact characterized by the preceding properties. Let $\Tc$ be a monad on $\Cs$ that is monoidal on $(\Cs,\diamond)$ and comonoidal on $(\Cs,\star)$. Suppose the unit map $\unita\diamond\unita\to\unita$ is bilinear and the interchange law is bilinear when postcomposed with a pair of bilinear maps, as in Lemma \ref{l:bimonoidal-duoidal}. 

\begin{lemma}\label{l:bimonoidal-duoidal2}
In the above situation, $\Tc$ is bimonoidal on $(\Cs,\diamond,\star)$.
\end{lemma}
\begin{proof}
As mentioned, we only need to verify the commutativity of \eqref{eq:bilaxass} and the second diagram in \eqref{eq:ibunit}. The latter is the bilinearity of $\unita\diamond\unita\to\unita$.
For the former one applies the hypothesis on $\zeta$ to the free algebras on objects $A$, $B$, $C$ and $D$, and the maps $\varphi_{A,B}$ and $\varphi_{C,D}$, which are bilinear by Lemma \ref{l:phi-multi}.
One precomposes with $\Tc(\eta_A\star\eta_B)\diamond\Tc(\eta_C\star\eta_D)$ to deduce the result.
\end{proof}

\subsection{Duoidal structure on the category of algebras}\label{ss:duoidal-alg}

Assume that universal $(\Tc,\varphi)$-multilinear maps in $(\Cs,\diamond)$ exist and that they are closed under composition. In other words, the multicategory $\Cs^{\Tc,\varphi}$ is representable by a monoidal structure $\diamond^\varphi$. The unit object is $\Tc\unit$. (Theorem \ref{t:monoidal} provides sufficient conditions for this, but we do not make the assumptions of the theorem. Lemma \ref{l:coeqmult} says that $\diamond^\varphi$ is defined in terms of Linton coequalizers, but we will not make use of this fact here.) Let
\[
u_{A_1,\dots,A_n}: A_1\diamond\dots\diamond A_n \to A_1\diamond^\varphi\dots\diamond^\varphi A_n
\]
denote the universal multilinear map in $(\Cs,\diamond)$ of a sequence of $\Tc$-algebras $A_1,\dots,A_n$. By assumption, diagrams of the form
\begin{equation}\label{eq:univ-comp}
\begin{gathered}
\xymatrix@C-20pt@R-10pt{
& A\diamond B\diamond C \ar[ld]_{u\diamond\id} \ar[dd]^-{u} \ar[rd]^{\id\diamond u} & \\
(A\diamond^\varphi B)\diamond C \ar[rd]_{u} & & A\diamond(B\diamond^\varphi C) \ar[ld]^{u}\\
& A\diamond^\varphi B\diamond^\varphi C &
}
\end{gathered}
\end{equation}
commute.


We proceed to turn $(\Cs,\diamond^\varphi,\Tc\unit,\star,\unita)$ into a duoidal category.
Let $A$, $B$, $C$ and $D$ be $\Tc$-algebras. By Lemma \ref{l:bimonoidal-duoidal}, the map
\[
(A \star B) \diamond (C \star D) \map{\zeta_{A,B,C,D}} (A \diamond C) \star (B \diamond D) \map{u_{A,C} \star u_{B,D}} (A \diamond^\varphi C) \star (B \diamond^\varphi D)
\]
is $(\Tc,\varphi)$-bilinear.
Universality yields an algebra morphism
\[
\zeta^\Tc_{A,B,C,D}: (A\star B) \diamond^\varphi (C\star D) \longrightarrow (A\diamond^\varphi C) \star (B \diamond^\varphi D)
\]
making the diagram
\begin{equation}\label{eq:zetaT}
\begin{gathered}
\xymatrix{
(A\star B) \diamond (C\star D) \ar[d]_-u \ar[r]^-\zeta & (A\diamond C) \star (B \diamond D)\ar[d]^-{u\star u} \\
(A\star B) \diamond^\varphi (C\star D) \ar@{-->}[r]_-{\zeta^\Tc} &(A\diamond^\varphi C) \star (B \diamond^\varphi D)
}
\end{gathered}
\end{equation}
commutative. This defines the interchange law $\zeta^\Tc$ for $\Cs^\Tc$. One of the unit maps is defined similarly, from the bilinearity of the corresponding unit map of $\Cs$:
\[
\xymatrix{
\unita\diamond\unita \ar[r] \ar[d]_{u_{\unita,\unita}} & \unita.\\
\unita\diamond^\varphi \unita \ar@{-->}[ru]
}
\]
The remaining unit maps are as follows.
\begin{gather*}
\Tc\unit \map{} \Tc(\unit\star\unit) \map{\psi_{\unit,\unit}} \Tc\unit\star\Tc\unit,\\
\Tc\unit \map{} \Tc\unita \map{\psi_0} \unita.
\end{gather*}
The unlabeled maps are built from the unit maps of $\Cs$.

\begin{theorem}\label{t:bimonoidal}
With the preceding structure, $(\Cs^\Tc, \diamond^\varphi, \Tc\unit, \star, \unita)$ is a duoidal category.
\end{theorem}
\begin{proof}
We verify one of the conditions for a duoidal category, namely axiom \eqref{eq:iass}.

We build the following diagram.
\[
\def\objectstyle{\scriptstyle}
\def\labelstyle{\scriptstyle}
\xymatrix@C+5pt@R-5pt{
(A\star X)\diamond(B\star Y)\diamond(C\star Z) \ar[ddd]_-{\zeta\diamond\id} \ar[dr]_{u\diamond\id} \ar[rr]^-{u} & & (A\star X)\diamond^\varphi(B\star Y)\diamond^\varphi(C\star Z)  \ar[ddd]^-{\zeta^\Tc\diamond^\varphi\id} \\
 & ((A\star X)\diamond^\varphi(B\star Y))\diamond(C\star Z)\ar[dd]_{\zeta^\Tc\diamond\id} \ar[ur]_{u} & \\
 & & \\
((A\diamond B)\star(X\diamond Y))\diamond(C\star Z) \ar[ddd]_-{\zeta} \ar[r]_{(u\star u)\diamond\id} & ((A\diamond^\varphi B)\star(X\diamond^\varphi Y))\diamond(C\star Z)\ar[dd]_{\zeta} \ar[r]_{u} & ((A\diamond^\varphi B)\star(X\diamond^\varphi Y))\diamond^\varphi(C\star Z)  \ar[ddd]^-{\zeta^\Tc} \\
 &  & \\
 &  ((A\diamond^\varphi B)\diamond C)\star((X\diamond^\varphi Y)\diamond Z)\ar[dr]^{u\star u}& \\
 (A\diamond B\diamond C)\star(X\diamond Y\diamond Z) \ar[rr]_-{u\star u}\ar[ur]^-{(u\diamond\id)\star(u\diamond\id)\ \ }   & & (A\diamond^\varphi B\diamond^\varphi C)\star(X\diamond^\varphi Y\diamond^\varphi Z)\\
}
\]
The triangles commute by \eqref{eq:univ-comp}, the quadrilaterals on the top left and bottom right corners by \eqref{eq:zetaT}, and the remaining two by naturality.

It follows that
\[
\zeta^\Tc\circ(\zeta^\Tc\diamond^\varphi\id)\circ u = (u\star u)\circ \zeta\circ(\zeta\diamond\id).
\]
A symmetric argument yields
\[
\zeta^\Tc\circ(\id\diamond^\varphi\zeta^\Tc)\circ u = (u\star u)\circ \zeta\circ(\id\diamond\zeta).
\]
Axiom \eqref{eq:iass} for $\zeta$ implies equality between the two previous expressions, and universality of $u$ allows us to conclude that
\[
\zeta^\Tc\circ(\zeta^\Tc\diamond^\varphi\id) = \zeta^\Tc\circ(\id\diamond^\varphi\zeta^\Tc). \qedhere
\]
\end{proof}

\subsection{Monoidal structure on the multicategory $\Cs^{\Tc,\varphi}$}\label{ss:monoidal-multi}
	
We retain the assumptions of Section \ref{ss:duoidal-alg}.
Recall the multicategory $\Cs^{\Tc,\varphi}$ from Section \ref{ss:multilinear}. It consists of $\Tc$-algebras and $(\Tc,\varphi)$-multilinear maps.  According to Lemma \ref{l:multisubcat}, there is a faithful morphism
\[
\Cs^{\Tc,\varphi} \to \Mcc(\Cs,\diamond), 
\]
where the latter is the multicategory represented by the monoidal structure $\diamond$ on $\Cs$. Proposition \ref{p:universal} states that $\Cs^{\Tc,\varphi}$ is represented by the monoidal structure $\diamond^\varphi$:
\begin{equation}\label{eq:multisubcat}
\Cs^{\Tc,\varphi} \cong \Mcc(\Cs^\Tc,\diamond^\varphi).
\end{equation}
A multilinear map $f$ is represented by the morphism of algebras $\tilde{f}$ defined by universality, as below.
\[
\xymatrix{
V_1\diamond\dots\diamond V_n \ar[d]_{u} \ar[r]^-{f} & W\\
V_1\diamond^\varphi\dots\diamond^\varphi V_n \ar[ru]_{\tilde{f}}
}
\]

The duoidal structure on $\Cs$ allow us to turn $\Mcc(\Cs,\diamond)$ into a \emph{monoidal multicategory}. We denote this structure by $\star_\zeta$. On objects we simply set
\[
X\star_\zeta Y = X\star Y.
\]
Given multimaps $f:X_1\diamond\dots\diamond X_n\to X$ and $g:Y_1\diamond\dots\diamond Y_n\to Y$, we let $f\star_\zeta g$ be the composite
\[
(X_1\star Y_1)\diamond\dots\diamond(X_n\star Y_n) \map{\zeta} 
(X_1\diamond\dots\diamond X_n)\star(Y_1\diamond\dots\diamond Y_n)\map{f\star g} X\star Y.
\]
The map on the left is the unique such obtained by iterating the interchange law $\zeta$ \cite[(6.19)]{AM}.

We may also consider the multicategory represented by $(\Cs^\Tc,\diamond^\varphi)$, and the monoidal structure $\star_{\zeta^\Tc}$ on it arising from the duoidal structure in Theorem \ref{t:bimonoidal}.

\begin{proposition}\label{p:monoidal-multi}
Let $(\Tc,\varphi,\psi)$ be a bimonoidal monad on $\Cs$. The multicategory
$\Cs^{\Tc,\varphi}$ inherits the monoidal structure $\star_\zeta$ of $\Mcc(\Cs,\diamond)$ under \eqref{eq:multisubcat}. Moreover, 
\[
(\Cs^{\Tc,\varphi},\star_\zeta) \cong (\Mcc(\Cs^\Tc,\diamond^\varphi), \star_{\zeta^\Tc})
\]
as monoidal multicategories.
\end{proposition}
\begin{proof}
The first statement is a generalization of Lemma \ref{l:bimonoidal-duoidal}, and may be established by similar arguments. For the second statement, we verify that tensor products of multimaps agree via the correspondence $f\leftrightarrow \tilde{f}$. Given algebras $A_i$ and $B_i$, we build the following diagram.
\[
\xymatrix@C-5pt{
(A_1\star B_1)\diamond\dots\diamond(A_n\star B_n) \ar[d]_{u} \ar[r]^-{\zeta} &
(A_1\diamond\dots\diamond A_n)\star(B_1\diamond\dots\diamond B_n) \ar[d]^{u\star u} \ar[r]^-{f\ast g} & C\star C\\
(A_1\star B_1)\diamond^\varphi\dots\diamond^\varphi(A_n\star B_n) \ar[r]_-{\zeta^\Tc} &
(A_1\diamond^\varphi\dots\diamond^\varphi A_n)\star(B_1\diamond^\varphi\dots\diamond^\varphi B_n) \ar[ru]_-{\tilde{f}\ast \tilde{g}}
}
\]
Its commutativity follows from \eqref{eq:zetaT}. We deduce that
\[
\widetilde{f\star_{\zeta} g} = \tilde{f}\star_{\zeta^\Tc} \tilde{g},
\]
as needed.
\end{proof}

\subsection{Example: linear monads}\label{ss:linear-bimonoidal-monad}
Let $(\Cs,\otimes)$ be a braided monoidal with braiding $\beta$. Let $H$ be a commutative bimonoid therein. The linear monad $H\otimes(-)$ on $\Cs$ is bimonoidal. The comonoidal structure is
\[
H\otimes X\otimes Y \map{\Delta\otimes\id} H\otimes H\otimes X\otimes Y \map{\id\otimes\beta\otimes\id}
H\otimes X\otimes H\otimes Y.
\]
Note this is an instance of the construction in Section \ref{ss:linear-monad}: $H$ is a bimonoid in the duoidal category associated to the braiding. The monoidal structure is as in Section \ref{ss:linear-monoidal-monad}.

Theorem \ref{t:bimonoidal} yields a duoidal structure on the category of $H$-modules
(under suitable assumptions such as those in Theorem \ref{t:monoidal}). The monoidal structures are $\diamond=\otimes_B$  (Example \ref{ss:linear-monoidal-monad}) and $\star=\otimes$ (Example \ref{ss:linear-monad}).

We illustrate with a special case of the dual to the preceding construction. Let $G$ be a monoid in the category of sets under cartesian product.
We may regard $G$ as a cocommutative bimonoid therein. (A fact available for all monoids in a cartesian category.) The linear comonad $G\times(-)$ is then bimonoidal and a duoidal structure on the category of $G$-graded sets ensues. It is the one described in Section \ref{ss:duoidal}.
For a related example, see Appendix \ref{s:species}.

This example may be generalized: starting from a $(2,1)$-monoid $H$ in a $3$-monoidal category $(\Cs,\diamond,\otimes,\star)$, the linear monad $H\otimes(-)$ is bimonoidal on the duoidal category $(\Cs,\diamond,\star)$. See Section \ref{ss:linear-higher} for more details.

\subsection{Example: pointing as a bimonoidal monad}\label{ss:pointing-bimon}

Let $(\Cs,\diamond,\unit,\star,\unita)$ be a duoidal category with interchange law $\zeta$. We name its unit maps as follows:
\[
\delta:\unit\to\unit\star\unit, \quad \zeta_0:\unit\to\unita, \quad \theta:\unita\diamond\unita\to\unita.
\]
We make the following main assumption: that the unit map $\zeta_0$ is an isomorphism. Duoidal categories with this property are called \emph{normal} \cite{BM2015,GL2016}. 

\begin{lemma}\label{l:duoidal-pointed}
In this situation, $\delta$ and $\theta$ are invertible.
\end{lemma}
\begin{proof}
Since $\unit$ is a comonoid with coproduct $\delta$ and counit $\zeta_0$, $\delta$ must be invertible with inverse $\zeta_0\star\id=\id\star\zeta_0$. The argument for $\theta$ is dual.
\end{proof}

As in Section \ref{ss:pointing-monad}, assume also that the unit object $\unit$ is terminal in $\Cs$, $\Cs$ admits coproducts and these are preserved by $X\diamond(-)$ and $(-)\diamond X$ for each object $X$. According to Lemma \ref{l:pointing-monoidal}, the pointing monad $\Tcp$ is then monoidal on $(\Cs,\diamond,\unit)$. It is also comonoidal on $(\Cs,\star,\unita)$; as discussed in Section \ref{ss:pointing-monad}, this relies on the fact that $\unit$ is a comonoid in $(\Cs,\star,\unita)$. 

\begin{lemma}\label{l:pointing-bimon}
In this situation, the pointing monad $\Tcp$ is bimonoidal on $(\Cs,\diamond,\unit,\star,\unita)$.
\end{lemma}
\begin{proof}
We verify the conditions of Lemma \ref{l:bimonoidal-duoidal2}, and for this we employ the criterion for $\Tcp$-bilinearity in \eqref{eq:smash-bilinear}.

Note first of all that $\zeta_0^{-1}$ must be the unique map $\unita\to\unit$. Then,
according to \eqref{eq:smash-bilinear}, bilinearity of $\theta$ is equivalent to the commutativity of the following diagrams.
\begin{equation*}
\begin{gathered}
\xymatrix@R-15pt{
\unita\diamond\unit \ar@{=}[d] \ar[r]^-{\id\diamond\zeta_0} & \unita\diamond\unita \ar[dd]^-{\theta}\\
\unita \ar[d]_{\zeta_0^{-1}} & \\
\unit \ar[r]_-{\zeta_0} & \unita
}
\end{gathered}
\qquad
\begin{gathered}
\xymatrix@R-15pt{
\unit\diamond\unita \ar@{=}[d] \ar[r]^-{\zeta_0\diamond\id} & \unita\diamond\unita \ar[dd]^-{\theta}\\
\unita \ar[d]_{\zeta_0^{-1}} & \\
\unit \ar[r]_-{\zeta_0} & \unita
}
\end{gathered}
\end{equation*}
This holds since $\unita$ is a monoid with multiplication $\theta$ and unit $\zeta_0$ (or by Lemma \ref{l:duoidal-pointed}).

Now let $A$, $B$, $C$, $D$, $X$ and $Y$ be pointed objects. We denote their points by $e$. Let $f:A\diamond C\to X$ and $g:B\diamond D\to Y$ be $\Tcp$-bilinear. We check that the composite $(f\star g)\circ\zeta$ is likewise. According to \eqref{eq:smash-bilinear}, bilinearity of $f$ and $g$ yields the commutative diagrams
\begin{equation*}
\begin{gathered}
\xymatrix@R-15pt@C+5pt{
 A\diamond\unit \ar@{=}[d] \ar[r]^-{\id_A\diamond e_C} & A\diamond C\ar[dd]^{f}\\
  A \ar[d] & \\
\unit \ar[r]_-{e_X} & X
}
\end{gathered}
\qquad
\begin{gathered}
\xymatrix@R-15pt@C+5pt{
 B\diamond\unit \ar@{=}[d] \ar[r]^-{\id_B\diamond e_D} & B\diamond D\ar[dd]^{g}\\
 B   \ar[d] & \\
\unit \ar[r]_-{e_Y} & Y,
}
\end{gathered}
\end{equation*}
(and two similar diagrams in which $\unit$ occurs on the left). The bilinearity of $(f\star g)\circ\zeta$  follows from the commutativity of the diagram below (and its left version).
\[
\xymatrix@R-10pt@C+25pt{
 (A\star B)\diamond\unit \ar@{=}[d] \ar[r]^-{\id\diamond\delta} & (A\star B)\diamond(\unit\star\unit) \ar[r]^-{\id\diamond(e_C\star e_D)} \ar[d]_{\zeta} & (A\star B)\diamond(C\star D) \ar[d]^{\zeta} \\
 A\star B \ar@{=}[r] \ar[dd] & (A\diamond\unit)\star(B\diamond\unit) \ar@{=}[d] \ar[r]_{(\id\diamond e_C)\star(\id\diamond e_D)} & (A\diamond C)\star(B\diamond D) \ar[dd]^{f\ast g}\\
  & A\star B \ar[d] & \\
 \unit \ar[r]_{\delta}  & \unit\star\unit \ar[r]_-{e_X\star e_Y} & X\star Y
}
\]
The square on the top left commutes by one of the axioms for a duoidal category. The square below it commutes since $\delta$ is invertible and $\unit$ is terminal.
\end{proof}

The discussion below makes use of some of the notions in Section \ref{s:double-comon}.

We now remove the assumption on terminality of $\unit$ in $\Cs$. We continue to assume that $\zeta_0$ is invertible and maintain the assumptions on existence and preservation of coproducts in $\Cs$. 

Assume in addition that products exist in $\Cs$. We proceed as in Section \ref{ss:smash}. The copointing comonad $\Tcc$ is defined on $\Cs$. Moreover, it is double monoidal on $(\Cs,\diamond,\star)$, in the dual sense to that of Section \ref{s:double-comon}. Indeed, the invertibility of $\zeta_0$ implies that $\unit$ is a double monoid in this duoidal category, and hence a  $(2,1)$-monoid in the $3$-monoidal category $(\Cs,\diamond,\star,\times)$. By (the dual of) Theorem \ref{t:doublecomonoidal}, the category of copointed objects $\Csc$ inherits the duoidal structure from $\Cs$. This category also inherits the coproducts from $\Cs$, and the copointed object $\unit$ is terminal in $\Csc$.
We now apply Lemma \ref{l:pointing-bimon} to the lifted monad $\Tcp$ on $(\Csc,\diamond,\unit,\star,\unita)$. We conclude that it is bimonoidal. Theorem \ref{t:bimonoidal} then yields a duoidal structure $(\diamond^\varphi,\star)$ on the category $\Cspc$ of bipointed objects (assuming representability as in Section \ref{ss:duoidal-alg}). Given bipointed objects $X$ and $Y$, $X\diamond^\varphi Y$ is their smash product, and $X\star Y$ is bipointed via
\[
\unit\map{\delta}\unit\star\unit \map{e_X\star e_Y} X\star Y \qand X\star Y \map{\epsilon_X\star\epsilon_Y} \unit\star\unit \map{\delta^{-1}}\unit.
\]

\section{Double monoidal monads}\label{s:double-mon}

There are three manners in which a monad may interact with the duoidal structure of the category on which it acts. One of these gives rise to the notion of bimonoidal monad studied
in Section~\ref{s:bimonoidal}. The other two occur when the
monad is either monoidal with respect to both monoidal structures, or comonoidal with respect to both structures.  In the latter case we speak of a double comonoidal monad. This notion is discussed in Section \ref{s:double-comon}.
In the former we speak of a double monoidal monad; we study this notion here. A central role is played by the notion of double multilinearity which we introduce.

Throughout, $(\Cs, \diamond, \unit, \star,\unita)$ denotes a duoidal category with interchange law $\zeta$.

\subsection{Double monoidal monads and duoidal structure on the category of algebras}\label{ss:double-mon}

A \emph{double monoidal monad} $(\Tc, \varphi, \gamma)$ on $\Cs$ consists of a monad $\Tc$ equipped with two structures $\varphi$ and $\gamma$ such that $(\Tc,\varphi)$ is monoidal on $(\Cs,\diamond,\unit)$, $(\Tc,\gamma)$ is monoidal on $(\Cs,\star,\unita)$, and the following diagrams commute.

\begin{equation}\label{eq:dlaxass}
\begin{gathered}
\xymatrix@C=-40pt@R=+20pt{
& (\Tc A \star \Tc B) \diamond (\Tc C \star \Tc D)
\ar[rd]^-{\zeta}
\ar[ld]_-{\gamma \diamond \gamma} & \\
{\Tc(A \star B) \diamond \Tc(C \star D)} \ar[d]_{\varphi}  & &
\bigl(\Tc A \diamond \Tc C\bigr) \star \bigl(\Tc B \diamond \Tc D\bigr)
\ar[d]^{\varphi \star \varphi}\\
{ \Tc\bigl((A \star B) \diamond (C \star D)\bigr)} \ar[rd]_-{\Tc \zeta} & & \Tc(A \diamond C) \star \Tc(B \diamond D)
\ar[dl]^-{\gamma} \\
& \Tc\bigl((A \diamond C) \star (B \diamond D)\bigr) 
}
\end{gathered}
\end{equation}


\begin{equation}\label{eq:dlunit}
\begin{gathered}
\xymatrix@R-0.5pc@C-1.2pc{
\unit \ar[d]_{}\ar[r]^-{\varphi_0} &
\Tc \unit \ar[r]^-{} &
\Tc(\unit \star \unit)\\
\unit \star \unit \ar[rr]_-{{\varphi}_0 \star {\varphi}_0} & &
\Tc \unit \star \Tc \unit \ar[u]_{\gamma}
}
\end{gathered}
\quad
\begin{gathered}
\xymatrix@R-0.5pc@C-1pc{
\unita \diamond \unita \ar[rr]^-{}  \ar[d]_-{{\gamma}_0 \diamond {\gamma}_0} & &
\unita \ar[d]^{{\gamma}_0}\\
\Tc\unita \diamond \Tc\unita \ar[r]_-{\varphi} & \Tc(\unita \diamond \unita) \ar[r]_-{}& \Tc \unita
}
\end{gathered}
\end{equation}

\begin{equation}\label{eq:dlinverse}
\begin{gathered}
\xymatrix@R-0.5pc@C-.5pc{
\unit \ar[r]_{} \ar[d]_-{{\varphi}_0} & \unita \ar[d]^{{\gamma}_0}\\
\Tc \unit \ar[r]_{} & \Tc \unita
}
\end{gathered}
\end{equation}
The unlabeled arrows are either unit maps of $\Cs$ or their images under $\Tc$.

Consider the $2$-category whose objects are duoidal categories, whose morphisms are double lax monoidal functors, and whose $2$-cells are morphisms of such functors \cite[Proposition 6.57]{AM}. A double monoidal monad is precisely a monad in this $2$-category.

We remark that \eqref{eq:dlinverse} and both diagrams in \eqref{eq:dlunit} are superfluous: the commutativity of these diagrams follows from the fact that $\eta$ is a morphism of monoidal functors for both $\varphi$ and $\gamma$, which includes the properties
\[
{\varphi}_0 = \eta_{\unit} \qand {\gamma}_0 = \eta_{\unita}.
\]

Apply the constructions of Section \ref{s:monoidal} to each of the monoidal monads $(\Tc,\varphi)$ and $(\Tc,\gamma)$. Under suitable assumptions  (Theorem \ref{t:monoidal}), we obtain two monoidal structures $\diamond^\varphi$ and $\star^\gamma$ on the category of $\Tc$-algebras.

\begin{theorem}\label{t:doublemonoidal}
Let $(\Tc, \varphi, \gamma)$ be a double monoidal monad on $\Cs$.
Assume that $\Cs$ admits reflexive coequalizers and that these are preserved by $\Tc$, $\diamond$ and $\star$.
Then $(\Cs^\Tc, \diamond^{\varphi}, \Tc\unit, \star^{\gamma}, \Tc\unita)$ is a duoidal category.
Moreover, the free algebra functor 
\[
\Fc^\Tc:(\Cs,\diamond, \star) \to (\Cs^\Tc,\diamond^{\varphi}, \star^{\gamma})
\] 
is double strong monoidal.
\end{theorem}

For the proof one adapts the results in Section~\ref{s:monoidal}. We sketch the main steps in the following sections. A crucial ingredient is the notion of \emph{double multilinearity} which we discuss next.

\subsection{Double multilinearity}\label{ss:double-multi}

 For clarity, we first spell out a special case.
Let $A$, $B$, $C$, $D$, $W$ be $\Tc$-algebras.
A map 
\[
f:(A \star B) \diamond (C \star D) \to W
\] 
in $\Cs$ is \emph{$(\Tc,\varphi,\gamma)$-bilinear} if the following diagram commutes.
\begin{equation}\label{eq:double-bilinear}
\begin{gathered}
\xymatrix@R-5pt@C+5pt{
(\Tc A \star \Tc B) \diamond (\Tc C \star \Tc D)\ar[ddd]_-{(a \star b) \diamond (b \star d)}   \ar[r]^-{\gamma \diamond \gamma} & \Tc(A \star B) \diamond \Tc (C \star D) \ar[d]^{\varphi}\\
 & \Tc \left ((A \star B) \diamond (C \star D) \right) \ar[d]^-{\Tc f}\\
 & \Tc W \ar[d]^-{w} \\
(A \star B) \diamond (C \star D) \ar[r]_-f & W
}
\end{gathered}
\end{equation}
To describe the general case, we set up some notation. Given a sequence
\[
\und{V} = (V_1,\dots,V_{m})
\]
of $\Tc$-algebras, let
\[
\star\und{V} = V_1\star\dots\star V_m, \quad \Tc\und{V} = (\Tc V_1,\dots,\Tc V_m),
\]
and
\[
\star\und{v} = v_1\star\dots\star v_m : \star \Tc\und{V} \to \star\und{V}
\]
be the $\star$-product of all the algebra structure maps in the sequence. Now let $\und{V^i}$ be a sequence of algebras, for each $i=1,\dots,n$, and let $W$ be another algebra. A map
\[
f: (\star\und{V^1})\diamond\dots\diamond (\star\und{V^n}) \to W
\]
in $\Cs$ is \emph{$(\Tc,\varphi,\gamma)$-multilinear} if the following diagram commutes.
\begin{equation}\label{eq:double-multilinear}
\begin{gathered}
\xymatrix@R-5pt@C+5pt{
(\star\Tc\und{V^1}) \diamond \dots\diamond (\star\Tc\und{V^n})
\ar[ddd]_-{\star\und{v}}   \ar[r]^-{\gamma \diamond\dots\diamond \gamma} & \Tc(\star\und{V^1}) \diamond \dots\diamond \Tc (\star\und{V^n}) \ar[d]^{\varphi}\\
 & \Tc \left ((\star\und{V^1}) \diamond \dots\diamond (\star\und{V^n}) \right) \ar[d]^-{\Tc f}\\
 & \Tc W \ar[d]^-{w} \\
(\star\und{V^1}) \diamond \dots\diamond (\star\und{V^n}) \ar[r]_-f & W
}
\end{gathered}
\end{equation}
Thus, \eqref{eq:double-bilinear} is the special case of \eqref{eq:double-multilinear} in which $n=2$, $\und{V^1}=(A,B)$~and~$\und{V^2}=(C,D)$. When all sequences are of length $1$, $f$ is simply a $(\Tc,\varphi)$-multilinear map in $(\Cs,\diamond)$. When there is only one sequence, $f$ is simply a $(\Tc,\gamma)$-multilinear map in $(\Cs,\star)$. When both occur, $f$ is a morphism of $\Tc$-algebras in $\Cs$. 

When $(\Tc,\varphi,\gamma)$ is understood, we speak of \emph{double multilinear} maps, or \emph{double bilinear} maps in the case of \eqref{eq:double-bilinear}.

The interchange law does not intervene in \eqref{eq:double-multilinear}. Double multilinearity may be considered for any ordered pair of monoidal structures on a category and on a monad. Though the first kind of maps play a more important role in our discussion, we consider both $(\Tc,\varphi,\gamma)$ and $(\Tc,\gamma,\varphi)$-multilinear maps. The interchange law turns maps of the second kind into maps of the first kind. We state the result in the case of bilinear maps, though the general statement and its proof are similar.

\begin{proposition}\label{p:dmultilinear}
Let $f: (A\diamond C)\star(B\diamond D)\to W$ be a $(\Tc,\gamma,\varphi)$-bilinear map. Then the composite
\[
(A\star B) \diamond (C \star D) \map{\zeta}
(A\diamond C)\star(B\diamond D) \map{f} W
\]
is $(\Tc,\varphi,\gamma)$-bilinear.
\end{proposition}
\begin{proof}
We build the following diagram.
\[
\xymatrix@C-5pt@R-5pt{
(A\star B)\diamond(C\star D)\ar[r]^\zeta &(A\diamond C)\star (B\diamond D) \ar[r]^-{f} & W\\
(\Tc A\star \Tc B)\diamond (\Tc C\star \Tc D)\ar[d]_{\gamma\diamond\gamma} \ar[u]^{(a\star b)\diamond(c\star d)}\ar[r]_-{\zeta} & (\Tc A\diamond \Tc C)\star (\Tc B\diamond \Tc D)   \ar[d]_{\varphi\star\varphi} \ar[u]^-{(a\diamond c)\star(b\diamond d)}& \\
\Tc(A\star B)\diamond \Tc  (C\star D) \ar[d]_{\varphi}& \Tc(A\diamond C)\star \Tc(B\diamond D) \ar[d]_{\gamma}\\
\Tc((A\star B)\diamond (C\star D))\ar[r]_-{\Tc \zeta}&
\Tc((A\diamond C)\star (B\diamond D))\ar[r]_-{\Tc f}& \Tc W\ar[uuu]_{w} \\
}
\]
The rectangle in the bottom left corner commutes by \eqref{eq:dlaxass} and the one on the right by $(\Tc,\gamma,\varphi)$-bilinearity of $f$.
\end{proof}

Conversely, if $\zeta$ satisfies the property of Proposition \ref{p:dmultilinear} for all such maps $f$, then axiom \eqref{eq:dlaxass} holds. Thus this property characterizes double monoidal monads. The proof is similar to that of Lemma \ref{l:bimonoidal-duoidal2}.

\subsection{Composition of double multilinear maps}\label{ss:comp-double-multi}

Again we first discuss a special case.

\begin{proposition}\label{p:dmonoidal-duoidal}
For any $(\Tc,\gamma)$-bilinear maps
\[
f: A\star B \to X \qand g: C\star D \to Y
\]
and any $(\Tc,\varphi)$-bilinear map
\[
h: X\diamond Y \to Z,
\]
the composite
\[
(A\star B) \diamond (C \star D) \map{f\diamond g}
 X\diamond Y \map{h} Z
\]
is $(\Tc,\varphi,\gamma)$-bilinear.
\end{proposition}
\begin{proof}
We build the following diagram.
\[
\xymatrix@R-5pt{
(A\star B)\diamond(C\star D) \ar[r]^-{f\star g} & X\diamond Y \ar[r]^-{h} & Z\\
(\Tc A\star \Tc B)\diamond (\Tc C\star \Tc D)\ar[d]_{\gamma\diamond\gamma} \ar[u]^{(a\star b)\diamond(c\star d)} & & \Tc Z\ar[u]_{z} \\
\Tc(A\star B)\diamond \Tc  (C\star D)  \ar[r]_-{\Tc f\diamond\Tc g} & \Tc X\diamond \Tc Y \ar[uu]_{x\diamond y} \ar[r]_-{\varphi} & \Tc (X\diamond Y) \ar[u]_-{\Tc h}
}
\]
The rectangles commute by bilinearity of $f$, $g$ and $h$.
\end{proof}

We turn to the details of the general case. There are two types of composition one may perform. In both cases, the result is $(\Tc,\varphi,\gamma)$-multilinear. 

\begin{enumerate}[(i)]
\item Given a sequence of $(\Tc,\varphi,\gamma)$-multilinear maps, we may form their $\diamond$-product and follow it with a $(\Tc,\varphi)$-multilinear map. 
\item We may start from a $(\Tc,\varphi,\gamma)$-multilinear map and precede it with a $\diamond$-product of maps, each of which is itself a $\star$-product of $(\Tc,\gamma)$-multilinear maps. 
\end{enumerate}

The composite in Proposition \ref{p:dmonoidal-duoidal} is a special case of both types of composition. 

Type (i) is as follows. Let 
\[
g:W_1\diamond\dots\diamond W_m\to X
\] 
be $(\Tc,\varphi)$-multilinear. For each $j=1,\dots,m$ and $i=1,\dots,n_j$, let $\und{V_j^i}$ be a sequence of algebras, and let
\[
f_j: (\star\und{V_j^1})\diamond\dots\diamond (\star\und{V_j^{n_j}}) \to W_j
\]
be a $(\Tc,\varphi,\gamma)$-multilinear map. The composite
\[
g\circ(f_1\diamond\dots\diamond f_m)
\]
is then $(\Tc,\varphi,\gamma)$-multilinear.

Type (ii) is as follows. For each $i=1,\dots,n$, let $\und{V^i}$ be a sequence of algebras and let
\[
g: (\star\und{V^1})\diamond\dots\diamond (\star\und{V^{n}}) \to W
\]
be $(\Tc,\varphi,\gamma)$-multilinear. Write $\und{V^i}=(V^i_1,\dots,V^i_{m_i})$. For each $i=1,\dots,n$ and $j=1,\dots,m_i$, let $\und{U^i_j}$ be a sequence of algebras, and let
\[
f^i_j: \star \und{U^i_j} \to V^i_j
\] 
be a $(\Tc,\gamma)$-multilinear map. Let $\und{U^i}$ denote the concatenation of the sequences 
$\und{U^i_1},\dots,\und{U^i_{m_i}}$ and let
\[
f^i : \star\und{U^i} \map{f^i_1\star\dots\star f^i_{m_i}} V^i_1\star\dots\star V^i_{m_i} = \star \und{V^i}.
\]
The composite
\[
(\star\und{U^1})\diamond\dots\diamond (\star\und{U^{n}}) \map{f^1\diamond\dots\diamond f^n}
(\star\und{V^1})\diamond\dots\diamond (\star\und{V^{n}}) \xrightarrow{g}W
\]
is then $(\Tc,\varphi,\gamma)$-multilinear.

\subsection{Universal double multilinear maps}\label{ss:univ-double-multi}

There is an evident notion of \emph{universality} for double multilinear maps. For example, a double bilinear map $\theta:(A \star B) \diamond (C \star D) \to U$ is \emph{universal} if for any bilinear map $f$ as above, there exists a unique morphism $\hat f$ in $\Cs^\Tc$ such that $\hat f \circ \theta=f$. When they exist, universal maps are unique up to isomorphism.

The following parallel pair in $\Cs^\Tc$ is reflexive:
\begin{equation}\label{eq:doublepair}
\xymatrix@C-8pt
{
\Tc \left ((\Tc A\star \Tc B)\diamond (\Tc C\star \Tc D)\right ) \ar@<0.5ex>[rrr]^-{\mu \circ \Tc\varphi \circ \Tc(\gamma\diamond\gamma)}
\ar@<-0.5ex>[rrr]_-{\Tc \left ((a\star b)\diamond (c\star d\right ))} &&&
\Tc \left ((A\star B)\diamond (C\star D)\right ).
}
\end{equation}
The common section is $\Tc(\eta_{A\star B} \diamond \eta_{C\star D})$.

In the same manner as $(\Tc, \varphi)$-multilinear maps are linked with Linton coequalizers \eqref{eq:pair}, $(\Tc,\varphi,\gamma)$-bilinear maps are linked with coequalizers of the pairs \eqref{eq:doublepair}.

\begin{lemma}\label{l:double-coeqmult}
Let $f\colon(A\star B)\diamond (C\star D) \to W$ be a map in $\Cs$. Then
\begin{enumerate}[(i)]
\item $f$ is $(\Tc,\varphi,\gamma)$-bilinear if and only if $\bar{f}$ coequalizes the pair \eqref{eq:doublepair}.
\item $f$ is universal $(\Tc,\varphi,\gamma)$-bilinear if and only if $\bar{f}$ is the coequalizer of the pair \eqref{eq:doublepair} in $\Cs^\Tc$.
\end{enumerate}
\end{lemma}

(As before, $\bar{f}:\Tc X\to V$ denotes the adjoint in $\Cs^\Tc$ of a map $f:X\to V$ in $\Cs$.)

The proof is similar to that of Lemma \ref{l:coeqmult}. The definition of the pair \eqref{eq:doublepair} can be extended to any number of sequences of algebras (instead of the two sequences $(A,B)$ and $(C,D)$), and Lemma \ref{l:double-coeqmult} holds in this generality.

\begin{proposition}\label{p:main-double}
Assume $\Cs$ admits reflexive coequalizers and these are preserved by $\Tc$, $\diamond$, and $\star$. Then universal multilinear maps are closed under compositions of either type (i) or (ii).
\end{proposition}

The proof is similar to that of Proposition \ref{p:main}. We outline the steps in the special case of a composite
\[
(A\star B)\diamond(C\star D) \map{u\diamond u} (A\star^\gamma B)\diamond(C\star^\gamma D) \map{v} (A\star^\gamma B)\diamond^\varphi(C\star^\gamma D)
\]
in which the maps $u$ are $(\Tc,\gamma)$-universal and the map $v$ is $(\Tc,\varphi)$-universal.
Let
\[
g:(A\star B)\diamond(C\star D)\to X
\]
be a $(\Tc,\varphi,\gamma)$-multilinear map. Let us abbreviate
\[
V = (A\star^\gamma B)\diamond(C\star^\gamma D) \qand W = (A\star^\gamma B)\diamond^\varphi(C\star^\gamma D).
\]
In the diagram
\[
\xymatrix@C+0pt{
\Tc \left((\Tc A\star \Tc B)\right)\diamond \Tc\left((\Tc C\star \Tc D)\right )  \ar[d]_{\varphi}  \ar@<0.5ex>[r]^-{} \ar@<-0.5ex>[r]_-{}  & \Tc(A\star B)\diamond \Tc(C\star D)  \ar[r]^-{\overline{u}\diamond\overline{u}} \ar[d]_{\varphi} &
 V \ar@{-->}[d]^{\tilde{g}} \ar[r]^-{v} & W \ar@{-->}[dl]^{\hat{g}} \\
\Tc \left ((\Tc A\star \Tc B)\diamond (\Tc C\star \Tc D)\right )  \ar@<0.5ex>[r]^-{} \ar@<-0.5ex>[r]_-{} & \Tc((A\star B)\diamond(C\star D))  \ar[r]_-{\bar{g}} & X
}
\]
the bottom pair is \eqref{eq:doublepair}, while the top pair is the $\diamond$-product of two pairs \eqref{eq:pair} for $\gamma$. Hence, $\bar{u}\diamond\bar{u}$ is the coequalizer of the top pair, while $\bar{g}$ coequalizes the bottom row (by Lemmas \ref{l:coeqmult} and \ref{l:double-coeqmult}). There is then a unique map $\tilde{g}$ as above.
By Lemma \ref{l:multilinear}, $\tilde{g}$ is $(\Tc,\varphi)$-multilinear. Then, universality of $v$ yields a (unique) morphism of algebras $\hat{g}$. Precomposing with $\eta_{A\star B}\diamond\eta_{C\star D}$ yields
\[
g = \hat{g}\circ v\circ(u\diamond u)
\]
which is the desired factorization. Uniqueness follows or one may appeal to Lemma \ref{l:univ-epi}.

We mention in passing that the result in Proposition \ref{p:main-double} holds for any pair of monoidal structures on a category $\Cs$, not necessarily linked by the duoidal axioms.

\subsection{The interchange law on algebras}\label{ss:law-dmonoidal}

Under the assumptions of Proposition~\ref{p:main-double},  universal bilinear maps exist both for $(\Tc,\varphi)$ and $(\Tc,\gamma)$. They are in particular universal $(\Tc,\varphi,\gamma)$-multilinear maps.
We use $u$ to denote a universal $(\Tc,\gamma)$-bilinear map and $v$ for a universal $(\Tc,\varphi)$-bilinear map.

We build the following diagram, given $\Tc$-algebras $A$, $B$, $C$ and $D$.
\begin{equation}\label{eq:zetaT-double}
\begin{gathered}
\xymatrix{
(A\star B) \diamond (C\star D) \ar[d]_{u\diamond u} \ar[r]^-\zeta & (A\diamond C) \star (B \diamond D)\ar[d]^-{v\star v} \\
(A\star^\gamma B) \diamond (C\star^\gamma D) \ar[d]_{v} & (A\diamond^\varphi C)\star (B \diamond^\varphi D) \ar[d]^{u}\\
(A\star^\gamma B) \diamond^\varphi (C\star^\gamma D) \ar@{-->}[r]_-{\zeta^\Tc} &(A\diamond^\varphi C) \star^\gamma (B \diamond^\varphi D)
}
\end{gathered}
\end{equation}
By Proposition \ref{p:main-double}, the vertical map on the left is universal $(\Tc,\varphi,\gamma)$-bilinear.
Combining Propositions \ref{p:dmultilinear} and \ref{p:dmonoidal-duoidal} we see that the composite around the top right corner is $(\Tc,\varphi,\gamma)$-bilinear.
Universality then defines a morphism $\zeta^\Tc$ in $\Cs^\Tc$ as indicated.
This is the interchange law of $\Cs^\Tc$. 

The unit maps of $\Cs^\Tc$ are obtained by applying $\Tc$ to the unit maps of $\Cs$.

In order to verify axiom \eqref{eq:iass} for $\Cs^\Tc$, we build the commutative diagram below. 
\[
\def\objectstyle{\scriptstyle}
\def\labelstyle{\scriptstyle}
\xymatrix@C-60pt{
(A\star X)\diamond(B\star Y)\diamond(C\star Z) \ar[ddd]_{\zeta\diamond\id} \ar[rd]_{u\diamond u\diamond\id} \ar[rr]^-{u\diamond u\diamond u} & & (A\star^\gamma X)\diamond(B\star^\gamma Y)\diamond(C\star^\gamma Z) \ar[rd]_{v\diamond\id} \ar[rr]^-{v} & & (A\star^\gamma X)\diamond^\varphi(B\star^\gamma Y)\diamond^\varphi(C\star^\gamma Z)\ar[ddddd]^{\zeta^\Tc\diamond^\varphi\id} \\
& ((A\star^\gamma X)\diamond(B\star^\gamma Y))\diamond(C\star Z) \ar[ru]_{\id\diamond u} \ar[rd]_{v\diamond\id} & & ((A\star^\gamma X)\diamond^\varphi(B\star^\gamma Y))\diamond(C\star^\gamma Z) \ar[ddd]^{\zeta^\Tc\diamond^\varphi\id} \ar[ru]_{v} & \\
& & ((A\star^\gamma X)\diamond^\varphi(B\star^\gamma Y))\diamond(C\star Z) \ar[ru]_{\id\diamond u} \ar[d]^{\zeta^\Tc\diamond\id} & & \\
((A\diamond B)\star(X\diamond Y))\diamond(C\star Z) \ar[ddd]_{\zeta} \ar[rd]|{(v\star v)\diamond\id} & &((A\diamond^\varphi B)\star^\gamma(X\diamond^\varphi Y))\diamond(C\star Z)   \ar[rd]^{\id\diamond u} & & \\
& ((A\diamond^\varphi B)\star(X\diamond^\varphi Y))\diamond(C\star Z) \ar[d]^{\zeta} \ar[ru]^{u\diamond\id} \ar[rr]_{u\diamond u} & &  ((A\diamond^\varphi B)\star^\gamma(X\diamond^\varphi Y))\diamond(C\star^\gamma Z) \ar[rd]^{v} & \\
& ((A\diamond^\varphi B)\diamond C)\star((X\diamond^\varphi Y)\diamond Z) \ar[rd]^{v\star v}  & & & ((A\diamond^\varphi B)\star^\gamma(X\diamond^\varphi Y))\diamond^\varphi(C\star^\gamma Z) \ar[d]^{\zeta^\Tc} \\
(A\diamond B\diamond C)\star(X\diamond Y\diamond Z)\ar[ru]|{(v\diamond\id)\star(v\diamond\id)} \ar[rr]_-{v\star v} & & (A\diamond^\varphi B\diamond^\varphi C)\star(X\diamond^\varphi Y\diamond^\varphi Z) \ar[rr]_-{u} & & (A\diamond^\varphi B\diamond^\varphi C)\star^\gamma(X\diamond^\varphi Y\diamond^\varphi Z)
}
\]
It follows that
\[
\zeta^\Tc\circ(\zeta^\Tc\diamond^\varphi\id)\circ v\circ(u\diamond u\diamond u) = u\circ (v\star v)\circ \zeta\circ(\zeta\diamond\id).
\]
A symmetric argument yields
\[
\zeta^\Tc\circ(\id\diamond^\varphi\zeta^\Tc)\circ v\circ(u\diamond u\diamond u) = u\circ (v\star v)\circ \zeta\circ(\id\diamond\zeta).
\]
Axiom \eqref{eq:iass} for $\zeta$ implies equality between the two previous expressions. Universality of  $v\circ(u\diamond u\diamond u)$ (Proposition \ref{p:main-double}) allow us to conclude that
\[
\zeta^\Tc\circ(\zeta^\Tc\diamond^\varphi\id) = \zeta^\Tc\circ(\id\diamond^\varphi\zeta^\Tc). \qedhere
\]

The remaining duoidal axioms may be verified similarly.

\subsection{The free algebra functor}\label{ss:free-dmonoidal}

Theorem \ref{t:monoidal} guarantees that the free algebra functor $\Fc^\Tc$ is strong monoidal for each of the monoidal structures $(\diamond^\varphi,\Tc\unit)$ and $(\star^\gamma,\Tc\unita)$. By Lemma~\ref{l:phi-multi}, 
\[
\varphi:\Tc X\diamond\Tc Y\to \Tc(X\diamond Y) \qand \gamma:\Tc X\star\Tc Y\to \Tc(X\star Y) 
\]
are universal $(\Tc,\varphi)$ and $(\Tc,\gamma)$-bilinear, respectively. Thus, we have isomorphisms $\bar{\varphi}$ and $\bar{\gamma}$ such that
\[
\xymatrix{
\Fc^\Tc X \diamond \Fc^\Tc Y \ar[r]^-{\varphi} \ar[d]_{u} & \Fc^\Tc(X\diamond Y)\\
\Fc^\Tc X \diamond^\varphi \Fc^\Tc Y \ar[ru]_{\bar{\varphi}}
}
\quad
\xymatrix{
\Fc^\Tc X \star \Fc^\Tc Y \ar[r]^-{\gamma} \ar[d]_{u} & \Fc^\Tc(X\star Y)\\
\Fc^\Tc X \star^\gamma \Fc^\Tc Y \ar[ru]_{\bar{\gamma}}
}
\]
commute. Let $A$, $B$, $C$ and $D$ be objects in $\Cs$. Consider diagram \eqref{eq:dlaxass}, which we redraw here as follows.
\[
\xymatrix@R-10pt{
(\Tc A \star \Tc B) \diamond (\Tc C \star \Tc D) \ar[r]^-{\zeta} \ar[d]_-{\gamma \diamond \gamma}  &
\bigl(\Tc A \diamond \Tc C\bigr) \star \bigl(\Tc B \diamond \Tc D\bigr)
\ar[d]^-{\varphi \star \varphi} \\
\Tc(A \star B) \diamond \Tc(C \star D) \ar[d]_-{\varphi}  & 
\Tc(A \diamond C) \star \Tc(B \diamond D)
\ar[d]^-{\gamma} \\
 \Tc\bigl((A \star B) \diamond (C \star D)\bigr) \ar[r]_-{\Fc^\Tc \zeta} &
\Tc\bigl((A \diamond C) \star (B \diamond D)\bigr) }
\]
Comparing with \eqref{eq:zetaT-double} we deduce the commutativity of the following diagram. 
\[
\xymatrix@R-10pt{
(\Fc^\Tc A \star^\gamma \Fc^\Tc B) \diamond^\varphi (\Fc^\Tc C \star^\gamma \Fc^\Tc D) \ar[r]^-{\zeta^\Tc} \ar[d]_-{\bar{\gamma} \diamond \bar{\gamma}}  &
\bigl(\Fc^\Tc A \diamond^\varphi \Fc^\Tc C\bigr) \star^\gamma \bigl(\Fc^\Tc B \diamond^\varphi \Fc^\Tc D\bigr)
\ar[d]^-{\bar{\varphi} \star \bar{\varphi}} \\
\Fc^\Tc(A \star B) \diamond^\varphi \Fc^\Tc(C \star D) \ar[d]_-{\bar{\varphi}}  & 
\Fc^\Tc(A \diamond C) \star^\gamma \Fc^\Tc(B \diamond D)
\ar[d]^-{\bar{\gamma}} \\
 \Fc^\Tc\bigl((A \star B) \diamond (C \star D)\bigr) \ar[r]_-{\Fc^\Tc \zeta} &
\Fc^\Tc\bigl((A \diamond C) \star (B \diamond D)\bigr) }
\]
This is precisely \cite[Axiom (6.35)]{AM} for the double strong monoidal functor $(\Fc^\Tc,\bar{\varphi},\bar{\gamma})$. 

This completes the proof of Theorem~\ref{t:doublemonoidal}.

\section{Double comonoidal monads}\label{s:double-comon}

We discuss this notion briefly. Let $(\Cs, \diamond, \unit, \star,\unita)$ be a duoidal category with interchange law $\zeta$.


A \emph{double comonoidal monad} $(\Tc, \chi, \psi)$ on $\Cs$ consists of a monad $\Tc$ equipped with two structures $\chi$ and $\psi$ such that $(\Tc,\chi)$ is comonoidal on $(\Cs,\diamond,\unit)$, $(\Tc,\psi)$ is comonoidal on $(\Cs,\star,\unita)$, and the following diagrams commute.

\begin{equation}\label{eq:dcaxass}
\begin{gathered}
\xymatrix@C=-40pt@R=+20pt{
&  \Tc\bigl((A \star B) \diamond (C \star D)\bigr) \ar[rd]^-{\Tc \zeta} \ar[ld]_-{\chi} & \\
{\Tc(A \star B) \diamond \Tc(C \star D)} \ar[d]_{\psi \diamond \psi}  & &  {\Tc\bigl((A \diamond C) \star (B \diamond D)\bigr)} \ar[d]^{\psi}\\
{ (\Tc A \star \Tc B) \diamond(\Tc C \star \Tc D)} \ar[rd]_-{\zeta} & & {\Tc(A \diamond C) \star \Tc(B \diamond D)}  \ar[dl]^-{\chi \star \chi} \\
& (\Tc A \diamond \Tc C) \star (\Tc B \diamond \Tc D)
}
\end{gathered}
\end{equation}


\begin{equation}\label{eq:dcunit}
\begin{gathered}
\xymatrix@R-0.5pc@C-1pc{
\Tc \unit\ar[d]_-{{\chi}_0} \ar[r]^-{} &
\Tc(\unit \star \unit) \ar[r]^{\psi}& \Tc\unit \star \Tc\unit   \ar[d]^-{{\chi}_0 \star {\chi}_0}\\
\unit \ar[rr]^-{}&&
\unit \star \unit
}
\quad
\xymatrix@R-0.5pc@C-1pc{
\Tc\unita\ar[d]_-{{\psi}_0} \ar[r]^-{} &
\Tc(\unita \diamond \unita)  \ar[r]^{\chi}  & \Tc\unita \diamond \Tc\unita \ar[d]^-{{\psi}_0 \diamond {\psi}_0}\\
\unita \ar[rr]_-{} && \unita \star \unita 
}
\end{gathered}
\end{equation}

\begin{equation}\label{eq:dcinverse}
\begin{gathered}
\xymatrix@R-0.5pc@C-.5pc{
\Tc\unit \ar[r]^{} \ar[d]_-{{\chi}_0} &
\Tc\unita \ar[d]^{{\psi}_0}\\
\unit \ar[r]_{} & \unita
}
\end{gathered}
\end{equation}

The unlabeled arrows are either unit maps of $\Cs$ or their images under $\Tc$.

Consider the $2$-category whose objects are duoidal categories, whose morphisms are double colax monoidal functors, and whose $2$-cells are morphisms of such functors \cite[Proposition 6.57]{AM}. A double comonoidal monad is precisely a monad in this $2$-category.

The following result interprets the axioms in terms of $\Tc$-algebras. Let $\Tc$ be a monad on $\Cs$ that is comonoidal on both $(\Cs,\diamond)$ and $(\Cs,\star)$. First note that each of the monoidal structures on $\Cs$ induces a corresponding structure on $\Cs^\Tc$, by Theorem \ref{t:comonoidal}. We do not distinguish notationally between the original and the induced structures. Given algebras $A$ and $B$, their tensor products are algebras via
\[
\Tc(A\diamond B)\map{\chi} \Tc A\diamond\Tc B\map{a\diamond b} A\diamond B
\qand
\Tc(A\star B)\map{\psi} \Tc A\star\Tc B\map{a\star b} A\star B.
\]
The unit objects $\unit$ and $\unita$ are algebras via $\chi_0$ and $\psi_0$, respectively.

\begin{lemma}\label{l:dcomonoidal-duoidal}
In the above situation, $\Tc$ is double comonoidal if and only all three unit maps of $\Cs$ are morphisms of $\Tc$-algebras, and for any $\Tc$-algebras $A$, $B$, $C$ and $D$, the interchange law $\zeta_{A,B,C,D}$ is a morphism of $\Tc$-algebras.
\end{lemma}
\begin{proof}
Each of the diagrams \eqref{eq:dcunit} and \eqref{eq:dcinverse} expresses precisely that one of the unit maps is a morphism of algebras. That $\zeta_{A,B,C,D}$ is a morphism of algebras (when the four objects are algebras) follows immediately from \eqref{eq:dcaxass}. Conversely, \eqref{eq:dcaxass} follows from the hypothesis on $\zeta$ applied to free algebras.
\end{proof}

The following is an immediate consequence.

\begin{theorem}\label{t:doublecomonoidal}
Let $(\Tc, \chi, \psi)$ be a double comonoidal monad on a duoidal category $\Cs$. Then $(\Cs^\Tc, \diamond, \unit, \star, \unita)$ is a duoidal category.
\end{theorem}

\section{Higher monoidal monads}\label{s:higher}

Bimonoidal, double monoidal, and double comonoidal monads constitute the second level in a hierarchy. For any $i$, $j$, and $n$ with $n=i+j$ and $i,j\geq 0$, one may consider
$(p,q)$-monoidal monads acting on an $n$-monoidal category.

\subsection{$3$-monoidal categories}\label{ss:3mon}

Let $\Cs$ be a category.
Three monoidal structures $(\diamond,\unit)$, $(\otimes,\unita)$ and $(\star,\unitaa)$ on $\Cs$, listed in order, turn $\Cs$ into a \emph{$3$-monoidal category} if the axioms in \cite[Definition 7.1.1]{AM} are satisfied. These include the requirements that each of 
$(\diamond,\unit,\otimes,\unita)$, 
$(\otimes,\unita,\star,\unitaa)$, and
$(\diamond,\unit,\star,\unitaa)$ is a duoidal structure on $\Cs$. We use $\zeta$ to denote the interchange law of any of them.
The axioms also require the commutativity of the following diagram.
\begin{equation}\label{eq:3braid}
\begin{gathered}
\def\objectstyle{\scriptstyle}
\def\labelstyle{\scriptstyle}
\xymatrix@R+1pc@C-1pc{
\left((A_1 \star B_1) \otimes (A_2 \star B_2)\right) \diamond
\left((C_1 \star D_1) \otimes (C_2 \star D_2)\right)
\ar[r]^-{\zeta \diamond \zeta}\ar[d]_{\zeta}
& \left((A_1 \otimes A_2) \star (B_1 \otimes B_2)\right)
\diamond \left((C_1 \otimes C_2) \star (D_1 \otimes D_2)\right)
\ar[d]^{\zeta}\\
\left((A_1 \star B_1) \diamond (C_1 \star D_1)\right) \otimes
\left((A_2 \star B_2) \diamond (C_2 \star D_2)\right)
\ar[d]_{\zeta \otimes \zeta}
& \left((A_1 \otimes A_2) \diamond (C_1 \otimes C_2)\right) \star
\left((B_1 \otimes B_2) \diamond (D_1 \otimes D_2)\right)
\ar[d]^{\zeta \star \zeta}\\
\left((A_1 \diamond C_1) \star (B_1 \diamond D_1)\right) \otimes
\left((A_2 \diamond C_2) \star (B_2 \diamond D_2)\right)
\ar[r]_-{\zeta}
& \left((A_1 \diamond C_1) \otimes (A_2 \diamond C_2)\right) \star
\left((B_1 \diamond D_1) \otimes (B_2 \diamond D_2)\right)
}
\end{gathered}
\end{equation}
The remaining axioms involve the unit maps of the three duoidal categories.

A symmetric monoidal structure on $\Cs$ gives rise to a $3$-monoidal structure in which $\diamond=\otimes=\star$ all coincide with the given structure. All three interchange laws also coincide and are constructed from the symmetry as in Section \ref{ss:duoidal}. See \cite[Section 7.3]{AM} for additional information.

$3$-monoidal structures may be constructed by combining a duoidal structure with a (co)cartesian structure on a category. If $(\Cs,\diamond,\star)$ is duoidal and $\times$ denotes the categorical product, then $(\Cs,\diamond,\star,\times)$ is $3$-monoidal. If $\amalg$ denotes the categorical coproduct, $(\Cs,\amalg,\diamond,\star)$ is $3$-monoidal. This is discussed (in part) in \cite[Section 7.3.1]{AM}. 

An example of a $3$-monoidal category in which none of the structures is cartesian or cocartesian is given in \cite[Proposition 8.69]{AM}.

\subsection{$n$-monoidal categories}\label{ss:n-mon}

Suppose given $n$ monoidal structures $(\otimes_i,\unit_i)$ on a category $\Cs$, and for each $i<j$, a duoidal structure on $(\Cs,\otimes_i,\unit_i,\otimes_j,\unit_j)$. We say that
\[
(\Cs,\otimes_1,\unit_1,\ldots,\otimes_n,\unit_n)
\]
is an \emph{$n$-monoidal category} if for each $i<j<k$,
\[
(\Cs,\otimes_i,\unit_i,\otimes_j,\unit_j,\otimes_k,\unit_k)
\]
is a $3$-monoidal category. The structure involves an interchange law and three unit maps for each pair $i<j$. See \cite[Section 7.6]{AM}.

If $(\Cs,\otimes_i,\unit_i)_{1\leq i\leq n}$ is $n$-monoidal and $\Cs$ admits (co)products, then
adding the (co)cartesian structure as the last (first) structure to the list produces an $(n+1)$-monoidal category.

\subsection{$(p,q)$-monoidal monads}\label{ss:pq-monad}

Let $(\Cs,\otimes_i,\unit_i)_{1\leq i\leq n}$ be an $n$-monoidal category. Fix $p$ and $q$ with $0\leq p,q\leq n$ and $p+q=n$. A \emph{$(p,q)$-monoidal monad} on $\Cs$ consists of a monad $\Tc$ equipped with $n$ structures $\gamma_i$ 
($1\leq i\leq n$) such that the first $p$ are monoidal, the last $q$ are comonoidal, and for each pair $i,j$ the monad $(\Tc, \gamma_i,\gamma_j)$ on the duoidal category $(\Cs,\otimes_i,\otimes_j)$ is:
\begin{itemize}
\item double monoidal if $i<j\leq p$,
\item bimonoidal if $i\leq p<j$,
\item double comonoidal if $p<i<j$.
\end{itemize}

In other words, a $(p,q)$-monoidal monad is a monad in the  $2$-category whose objects are $n$-monoidal functors and whose morphisms are $(p,q)$-monoidal functors \cite[Section 7.8]{AM}.

\subsection{Higher monoidal structure on the category of algebras}\label{ss:higher-alg}

Apply the constructions of Section \ref{s:monoidal} to each of the monoidal structures $\gamma_i$ on $\Tc$, $1\leq i\leq p$.  Under suitable assumptions  (Theorem \ref{t:monoidal}), we obtain monoidal structures $\otimes^{\gamma_i}$ on the category of $\Tc$-algebras. On the other hand, according to Theorem \ref{t:comonoidal}, for each $p<i$ the comonoidal structure $\gamma_i$ allow us to lift the monoidal structure $\otimes_i$ from $\Cs$ to $\Cs^\Tc$. 

\begin{theorem}\label{t:nmonoidal}
Assume that $\Cs$ admits reflexive coequalizers and that these are preserved by $\Tc$ and $\otimes_i$ for each $i\leq p$.  Then 
\[
(\Cs^\Tc, \otimes^{\gamma_1}, \dots,\otimes^{\gamma_p}, \otimes_{p+1}, \dots, \otimes_n)
\] 
is an $n$-monoidal category.
\end{theorem}

We outline the main steps in the proof.

The interchange law between the $i$-th and $j$-th monoidal structures on $\Cs$ is afforded by either Theorem \ref{t:doublemonoidal} if $i<j\leq p$, Theorem \ref{t:bimonoidal} if $i\leq p< j$, or Theorem \ref{t:doublecomonoidal} if $p<i<j$. One then has to verify the $3$-monoidal axioms for each triple $i<j<k$. There are four possible cases.

\begin{enumerate}[(i)]
\item $p<i<j<k$. In this case, the axioms follow from the fact that $\Uc^\Tc$ is strong monoidal with respect to all three monoidal structures.
\item $i\leq p< j<k$. The axioms follow by arguments involving universal $(\Tc,\gamma_i)$-multilinear maps similar to that in the proof of Theorem~\ref{t:bimonoidal}.
\item $i<j\leq p<k$. 
The axioms now involve $(\Tc,\gamma_i,\gamma_j)$-bilinear maps and are similar to those in the proof of Theorem~\ref{t:doublemonoidal}.
\item $i<j<k\leq p$. We provide more detail in this, the case that requires more
  attention.
\end{enumerate}

Let $u_i$, $u_j$ and $u_k$ denote universal bilinear maps for each of the monoidal structures.
Let $U_{ijk}$ denote the composite below.
\[
\xymatrix{
\left((A_1 \otimes_k B_1) \otimes_j (A_2 \otimes_k B_2)\right) \otimes_i
\left((C_1 \otimes_k D_1) \otimes_j (C_2 \otimes_k D_2)\right) \ar[d]^{(u_k \otimes_j u_k)\otimes_i (u_k \otimes_j u_k)}\\
\left((A_1 \otimes^{\gamma_k} B_1) \otimes_j (A_2 \otimes^{\gamma_k} B_2)\right) \otimes_i
\left((C_1 \otimes^{\gamma_k} D_1) \otimes_j (C_2 \otimes^{\gamma_k} D_2)\right) \ar[d]^{u_j \otimes_i u_j}\\
\left((A_1 \otimes^{\gamma_k} B_1) \otimes^{\gamma_j} (A_2 \otimes^{\gamma_k} B_2)\right) \otimes_i
\left((C_1 \otimes^{\gamma_k} D_1) \otimes^{\gamma_j} (C_2 \otimes^{\gamma_k} D_2)\right)
\ar[d]^{u_i}\\
\left((A_1 \otimes^{\gamma_k} B_1) \otimes^{\gamma_j} (A_2 \otimes^{\gamma_k} B_2)\right) \otimes^{\gamma_i}
\left((C_1 \otimes^{\gamma_k} D_1) \otimes^{\gamma_j} (C_2 \otimes^{\gamma_k} D_2)\right)
}
\]
The map $U_{kji}$ is defined in the same manner.

Let $\zeta_{ij}$, $\zeta_{ik}$, $\zeta_{jk}$ be the interchange laws for $(\Cs,\otimes_i,\otimes_j,\otimes_k)$, and $\zeta_{ij}^\Tc$, $\zeta_{ik}^\Tc$, $\zeta_{jk}^\Tc$ those for $(\Cs^\Tc,\otimes^{\gamma_i},\otimes^{\gamma_j},\otimes^{\gamma_k})$. Arguing as in the above proofs, and employing axiom \eqref{eq:3braid} for $\Cs$, one obtains
\begin{multline*}
\zeta^\Tc_{jk}\circ(\zeta^\Tc_{ik}\otimes_j\zeta^\Tc_{ik})\circ\zeta^\Tc_{ij} \circ U_{ijk}
= U_{kji} \circ \zeta_{jk}\circ(\zeta_{ik}\otimes_j\zeta_{ik})\circ\zeta_{ij} \\
=U_{kji} \circ(\zeta_{ij}\otimes_k\zeta_{ij})\circ\zeta_{ik}\circ(\zeta_{jk}\otimes_i\zeta_{jk}) 
= (\zeta^\Tc_{ij}\otimes^{\gamma_k}\zeta^\Tc_{ij})\circ\zeta^\Tc_{ik}\circ(\zeta^\Tc_{jk}\otimes^{\gamma_i}\zeta^\Tc_{jk})\circ U_{ijk}.
\end{multline*}
Thus, in order to to obtain axiom \eqref{eq:3braid} for $\Cs^\Tc$, it suffices to check that $U_{ijk}$ is epimorphic with respect to morphisms of $\Tc$-algebras.
This holds since $u_i$ satisfies this property (by universality) and since $(u_j \otimes_i u_j) \circ((u_k \otimes_j u_k)\otimes_i (u_k \otimes_j u_k))$ is epimorphic with respect to $(\Tc,\gamma_i)$-bilinear maps. This is the
case by Lemma \ref{l:univ-epi}, since $u_j\circ(u_k\otimes_j u_k)$ is a universal
$(\Tc,\gamma_j,\gamma_k)$-multilinear map (Proposition \ref{p:main-double}), hence
a reflexive $\Tc$-coequalizer (Lemma \ref{l:double-coeqmult}).

\subsection{$(p,q)$-monoids}\label{ss:pq-monoid}

Let $(\Cs,\otimes_i,\unit_i)_{1\leq i\leq n}$ be an $n$-monoidal category. Fix $p$ and $q$ with $0\leq p,q\leq n$ and $p+q=n$. A \emph{$(p,q)$-monoid} on $\Cs$ consists of an object $H$ equipped with a monoid structure in $(\Cs,\otimes_i)$ for each $i\leq p$ and a comonoid structure
in $(\Cs,\otimes_i)$ for each $i>p$, such that for each $i<j$, the $i$-th and $j$-th structures turn $M$ into a
\begin{itemize}
\item double monoid if $i<j\leq p$,
\item bimonoid if $i\leq p<j$,
\item double comonoid if $p<i<j$,
\end{itemize}
in the duoidal category $(\Cs,\otimes_i,\otimes_j)$. See \cite[Section 7.7]{AM}.

\subsection{Linear monads on higher monoidal categories}\label{ss:linear-higher}

Let $n=p+q$ with $p,q\geq 0$.
Suppose now that $H$ is a $(p+1,q)$-monoid in an $(n+1)$-monoidal category $(\Cs,\otimes_i)_{1\leq i\leq n+1}$.
We employ the last of the $p+1$ monoid structures to define a monad $\Tc$ on $\Cs$ by
\[
\Tc X = H\otimes_{p+1} X.
\]
We employ the remaining $n$ structures to endow $\Tc$ with $p$ monoidal and $q$ comonoidal structures: for $i<p+1<j$,
\begin{gather*}
(H\otimes_{p+1} X)\otimes_i (H\otimes_{p+1} Y) \map{\zeta} (H\otimes_i H) \otimes_{p+1} (X\otimes_i Y) \to H\otimes_{p+1}(X\otimes_i Y),\\
H\otimes_{p+1}(X\otimes_j Y) \to  (H\otimes_j H) \otimes_{p+1} (X\otimes_i Y)  \map{\zeta} (H\otimes_{p+1} X)\otimes_j (H\otimes_{p+1} Y).
\end{gather*}

\begin{proposition}\label{p:linear-higher}
With this structure, $\Tc$ is a $(p,q)$-monoidal monad on the $n$-monoidal category $(\Cs,\otimes_i)_{i\neq p+1}$.
\end{proposition}
\begin{proof}
Given an object $A$, let $\Fc_A:\Is\to\Cs$ denote the functor that sends the unique object of the one-arrow category $\Is$ to $A$. Then $A$ is a $(p,q)$-monoid if and only $\Fc_A$ is $(p,q)$-monoidal \cite[Proposition 7.40]{AM}. On the other hand, the definition of higher monoidal category implies that $\otimes_{p+1}:\Cs\times\Cs\to\Cs$ is $(p,q)$-monoidal. It follows that, if $A$ is such a monoid, the composite
\[
\Cs = \Is\times\Cs \map{\Fc_A\times\id} \Cs\times\Cs \map{\otimes_{p+1}} \Cs
\]
is $(p,q)$-monoidal. We obtain a functor
\[
A\mapsto A\otimes_{p+1}(-)
\]
from the category of $(p,q)$-monoids in $(\Cs,\otimes_i)_{i\neq p+1}$ to the category of $(p,q)$-monoidal endofunctors of this $n$-monoidal category. The former is monoidal under $\otimes_{p+1}$ and the latter under functor composition, and the above functor is strong monoidal for these structures. The result follows since a $(p+1,q)$-monoid is a monoid in the former and a $(p,q)$-monoidal monad is a monoid in the latter. 
\end{proof} 

Such $(p,q)$-monoidal monads are called \emph{linear}.

Algebras over $\Tc$ are modules over $H$ in $(\Cs,\otimes_{p+1})$.
Under the assumptions of Theorem \ref{t:nmonoidal}, we obtain an $n$-monoidal structure on this category.

The dual construction inputs a $(p,1+q)$-monoid $H$ in an $(n+1)$-monoidal category $\Cs$ and produces the $(p,q)$-comonad $H\otimes_{p+1}(-)$.

We list the possibilities for a linear monad and a linear comonad on a $3$-monoidal category
$(\Cs,\diamond,\otimes,\star)$. Note the former necessitates at least one monoid structure on $H$, the latter at least one comonoid structure.

\begin{table}[!h]
\begin{tabular}{@{}c|c|c|c@{}}
\hline
monoid type & linear monad & acted category & monoidality\\
\hline
$(3,0)$ & $H\star(-$) & $(\Cs,\diamond,\otimes)$ &  double monoidal\\
$(2,1)$ & $H\otimes(-)$ & $(\Cs,\diamond,\star)$ & bimonoidal\\
$(1,2)$ & $H\diamond(-)$ & $(\Cs,\otimes,\star)$ & double comonoidal\\
\hline
\end{tabular}
\end{table}

\begin{table}[!h]
\begin{tabular}{@{}c|c|c|c@{}}
\hline
monoid type & linear comonad & acted category & monoidality\\
\hline
$(0,3)$ & $H\diamond(-)$ & $(\Cs,\otimes,\star)$ & double comonoidal\\
$(1,2)$ & $H\otimes(-)$ & $(\Cs,\diamond,\star)$ & bimonoidal\\
$(2,1)$ & $H\star(-)$ & $(\Cs,\diamond,\otimes)$ & double monoidal\\
\hline
\end{tabular}
\end{table}

We encountered an example of a double monoidal linear comonad in Section \ref{ss:pointing-bimon} (the copointing comonad). In Appendix \ref{s:species} we describe certain bimonoidal comonads on the category of species which though not quite linear, are constructed from similar data on a given species.

\appendix

\section{Comonads on the category of species}\label{s:species}

We discuss a construction of bimonoidal comonads on the category of species.
It shares some features with the construction of Section \ref{ss:linear-bimonoidal-monad},
but it is set apart by the fact that the present examples are non-linear. More precisely, while the underlying monoidal comonads are linear, their comonoidal structure is not. We expand on this below.

Let $\Ss$ denote Joyal's category of \emph{set species}. We employ the terminology of \cite[Chapter 8]{AM} and \cite[Chapter 4]{AM13}. Warning: in the former reference, $\Ss$ denotes the category of \emph{vector species}. The distinction between the two categories plays a role below. 

A (set) species is a functor
\[
\Fset\to\Set
\]
from the groupoid of finite sets and bijections to the category of sets. 
A species $\rP$ specifies a set $\rP[I]$ for each finite set $I$. 
A morphism of species is a natural transformation $f:\rP\to\rQ$, and as such has components
\[
f_I:\rP[I]\to\rQ[I],
\]
one for each finite set $I$. We will not dwell on the role played by the morphisms in $\Fset$.

The category of species is symmetric monoidal under \emph{Cauchy product}: the species $\rP\bdot\rQ$ is given on a finite set $I$ by
\[
(\rP\bdot\rQ)[I] = \coprod_{I=S\sqcup T} \rP[S]\times\rQ[T].
\]
The disjoint union is taken over all ordered pairs $(S,T)$ of subsets of $I$ such that $I=S\cup T$ and $S\cap T=\emptyset$. The unit object is the species $\rU$ characteristic of the empty set. The symmetry simply exchanges the factors in the product.

Let $\rB$ be a \emph{set-theoretic bimonoid} in $\Ss$.
Thus, for each finite set $I$ and each decomposition $I=S\sqcup T$, there are maps
\[
\mu_{S,T} : \rB[S]\times\rB[T] \to \rB[I] \qqand \Delta_{S,T}:\rB[I]\to\rB[S]\times\rB[T].
\]
There is also a distinguished element $1\in\rB[\emptyset]$. The structure is subject to the axioms described in \cite[Section 4.2]{AM13}. While the maps $\mu$ and the element $1$ turn $\rB$ into a monoid in the monoidal category $(\Ss,\bdot,\rU)$, the maps $\Delta$ do not endow it with the structure of a comonoid therein. Note that if there is a morphism $\rP\to\rU$ in $\Ss$, then $\rP[I]=\emptyset$ on every nonempty set $I$. Certainly we do not require this from $\rB$.

There are many natural examples of set-theoretic bimonoids. 
The species $\rL$ of \emph{linear orders} is one: given linear orders $\ell_1$ on $S$ and $\ell_2$ on $T$, their concatenation defines a linear order $\mu_{S,T}(\ell_1,\ell_2)$ on $I$. Given a linear order $\ell$ on $I$, their restrictions define $\Delta_{S,T}(\ell)=(\ell|_S,\ell|_T)$.
Another example is furnished by the species of \emph{matroids}: direct sum of matroids defines the multiplication, while the matroid operations of restriction and contraction define the comultiplication.
See \cite[Chapters 12 and 13]{AM} and \cite[Section 9]{AM13} for more examples.

We employ a second monoidal structure on $\Ss$, the \emph{Hadamard product}: the species $\rP\bdot\rQ$ is given on a finite set $I$ by
\[
(\rP\times\rQ)[I] =  \rP[I]\times\rQ[I].
\]
Note this is simply the categorical product on $\Ss$. Therefore,
together with the Cauchy product, it turns the category of species into a duoidal category $(\Ss,\bdot,\times)$. Moreover, a bimonoid in this duoidal category is the same as a monoid in $(\Ss,\bdot)$.

Consider in addition the duoidal category $(\Ss,\bdot,\bdot)$ associated to the symmetric monoidal $(\Ss,\bdot)$ (Section \ref{ss:duoidal}).
We construct a bimonoidal comonad $\Tc$ on this duoidal category out of the above structure on $\rB$. The endofunctor is defined by
\[
\Tc \rP = \rB\times\rP.
\]
The comultiplication $\Tc\rP\to\Tc^2\rP$ and counit $\Tc\rP\to\rP$ of the comonad have components as follows:
\begin{align*}
\rB[I]\times\rP[I] &\to \rB[I]\times \rB[I]\times\rP[I] & & & \rB[I]\times\rP[I] &\to \rP[I]\\
(a,x) & \mapsto (a,a,x) & & & (a,x) &\mapsto x.
\end{align*}
The monoidal structure $\varphi_{\rP,\rQ}:\Tc \rP \bdot \Tc \rQ \to \Tc(\rP\bdot \rQ)$ has components
\begin{align*}
\rB[S]\times\rP[S]\times\rB[T]\times\rB[T] &\to \rB[I]\times\rP[S]\times\rQ[T] \\
(a,x,b,y) &\mapsto (\mu_{S,T}(a, b), x,y).
\end{align*}
The map $\varphi_0:\rU\to\Tc\rU$ has only one non-trivial component: 
\begin{align*}
\rU[\emptyset] & \to \rB[\emptyset]\times\rU[\emptyset]\\
\ast & \mapsto (1,\ast).
\end{align*}
The comonoidal structure $\psi_{\rP,\rQ}:\Tc(\rP\bdot \rQ)\to \Tc \rP \bdot \Tc \rQ$ has components
\begin{align*}
 \rB[I]\times\rP[S]\times\rQ[T] & \to \rB[S]\times\rP[S]\times\rB[T]\times\rQ[T] \\
(a,x,y) & \mapsto (a_1, x,a_2,y),
\end{align*}
where $\Delta_{S,T}(a)=(a_1,a_2)$. The map $\psi_0:\Tc\rU\to\rU$ has only one non-trivial component:
\begin{align*}
 \rB[\emptyset]\times\rU[\emptyset] & \to \rU[\emptyset] \\
(x,\ast)  & \mapsto \ast.
\end{align*}

These definitions turn $\Tc$ into a bimonoidal comonad on $(\Ss,\bdot,\bdot)$. The axioms
follow from the axioms for the set-theoretic bimonoid $\rB$. 

We now describe in explicit terms the duoidal structure on the category of $\Tc$-comodules afforded by Theorem \ref{t:bimonoidal}. 

To this end, first recall the \emph{groupoid of elements} $\el(\rB)$ associated to the species $\rB$. The objects are pairs $[I,a]$ where $I$ is a finite set and $a$ is an element of $\rB[I]$. A morphism $[I,a]\to[J,b]$ in $\rB$ is a morphism $\sigma:I\to J$ in $\Fset$ such that $\rB[\sigma](a)=b$.
A \emph{$\rB$-species} is a functor
\[
\el(\rB)\to\Set.
\]
Let $\Ss_\rB$ denote the category of $\rB$-species.

Note that the category of $\Tc$-comodules is the slice category of $\Ss$ under the species $\rB$. We may therefore identify it with $\Ss_\rB$. Concretely: a $\Tc$-coalgebra is a species $\rP$ equipped with a morphism $f:\rP\to\rB$ in $\Ss$. To this one associates the $\rB$-species
\[
\el(\rB)\to\Set, \quad [I,a]\mapsto \{x\in\rP[I] \mid f_I(x)=a\}.
\]

The resulting duoidal structure $(\diamond,\star)$ on the category $\Ss_\rB$ is as follows. Given two $\rB$-species $\rP$ and $\rQ$, we have 
\[
(\rP\diamond\rQ)[I,a] = \coprod_{\substack{I=S\sqcup T\\b\in\rB[S],\,c\in\rB[T]\\\mu_{S,T}(b,c)=a} }\rP[S,b]\times\rQ[T,c]
\]
and
\[
(\rP\star\rQ)[I,a] = \coprod_{I=S\sqcup T} \rP[S,a_1]\times\rQ[T,a_2],
\]
where $\Delta_{S,T}(a)=(a_1,a_2)$. The respective unit objects $\rI$ and $\rJ$ have
\[
\rI[I,a]=
\begin{cases}
\{\ast\} & \text{ if $I=\emptyset$ and $a=1$,} \\
\emptyset         & \text{ if not,}
\end{cases}
\qquad
\rJ[I,a] =
\begin{cases}
\{a\} & \text{ if $I=\emptyset$} \\
\emptyset         & \text{ if not.}
\end{cases}
\]

We close with some remarks regarding the non-linearity of the bimonoidal comonad $\Tc$.
The structure of monoidal comonad on $\Tc$ is in fact linear: $\rB$ is a bimonoid in $(\Ss,\bdot,\times)$ (being a monoid in $(\Ss,\bdot)$) and $\Tc$ is the associated linear monoidal comonad (in the dual sense to that of Section \ref{ss:linear-monad}).
On the other hand, the comonoidal structure is not linear, since as explained, $\rB$ is not a comonoid in $(\Ss,\bdot)$ (unless the species $\rB$ is concentrated on the empty set).

Note also that for $\Tc$ to be a linear bimonoidal comonad (in the dual sense to that of Section \ref{ss:linear-bimonoidal-monad}), $\rB$ would have to be a $(2,1)$-monoid in a $3$-monoidal category of the form $(\Ss,\bdot,\times,\bdot)$ (see Section \ref{ss:linear-higher}).
However, not even a duoidal structure of the form $(\Ss,\times,\bdot)$ exists: there is no natural transformation of the appropriate form for an interchange law.

The situation changes if we work with the category $\Ss_\Kb$ of \emph{vector species} instead: we do have a $3$-monoidal category $(\Ss_\Kb,\bdot,\times,\bdot)$ \cite[Proposition 8.69]{AM}, a set-theoretic bimonoid as above gives rise by linearization to a  $(2,1)$-monoid therein, and we have an associated bimonoidal comonad. The resulting duoidal structure on $\Ss_\Kb$ admits a description similar to the above.

\bibliographystyle{abbrv}
\bibliography{monads}

\newcommand{\noopsort}[1]{} \newcommand{\singleletter}[1]{#1} \def\cprime{$'$}
\begin{thebibliography}{10}

\bibitem{ARV2011}
J.~Ad{\'a}mek, J.~Rosick{\'y}, and E.~M. Vitale.
\newblock {\em Algebraic theories}, volume 184 of {\em Cambridge Tracts in
  Mathematics}.
\newblock Cambridge University Press, Cambridge, 2011.
\newblock A categorical introduction to general algebra, With a foreword by F.
  W. Lawvere.

\bibitem{AC}
M.~Aguiar and S.~U. Chase.
\newblock Generalized {H}opf modules for bimonads.
\newblock {\em Theory Appl. Categ.}, 27:263--326, 2012.

\bibitem{AM}
M.~Aguiar and S.~Mahajan.
\newblock {\em Monoidal functors, species and {H}opf algebras}, volume~29 of
  {\em CRM Monograph Series}.
\newblock American Mathematical Society, Providence, RI, 2010.
\newblock With forewords by Kenneth Brown and Stephen Chase and Andr{\'e}
  Joyal.

\bibitem{AM13}
M.~Aguiar and S.~Mahajan.
\newblock Hopf monoids in the category of species.
\newblock In {\em Hopf algebras and tensor categories}, volume 585 of {\em
  Contemp. Math.}, pages 17--124. Amer. Math. Soc., Providence, RI, 2013.

\bibitem{AJ}
M.~Anel and A.~Joyal.
\newblock Sweedler {T}heory for (co)algebras and the bar-cobar constructions,
  available{\noopsort{2013}} at \url{arXiv:1309.6952}.

\bibitem{BFSV}
C.~Balteanu, Z.~Fiedorowicz, R.~Schw{\"a}nzl, and R.~Vogt.
\newblock Iterated monoidal categories.
\newblock {\em Adv. Math.}, 176(2):277--349, 2003.

\bibitem{BW}
M.~Barr and C.~Wells.
\newblock {\em Toposes, triples and theories}, volume 278 of {\em Grundlehren
  der Mathematischen Wissenschaften [Fundamental Principles of Mathematical
  Sciences]}.
\newblock Springer-Verlag, New York, 1985.

\bibitem{BM2015}
M.~Batanin and M.~Markl.
\newblock Operadic categories and duoidal {D}eligne's conjecture.
\newblock {\em Adv. Math.}, 285:1630--1687, 2015.

\bibitem{BCZ2013}
G.~B{\"o}hm, Y.~Chen, and L.~Zhang.
\newblock On {H}opf monoids in duoidal categories.
\newblock {\em J. Algebra}, 394:139--172, 2013.

\bibitem{BGL2014}
G.~B{\"o}hm, J.~G{\'o}mez-Torrecillas, and E.~L{\'o}pez-Centella.
\newblock On the category of weak bialgebras.
\newblock {\em J. Algebra}, 399:801--844, 2014.

\bibitem{BS2013}
T.~Booker and R.~Street.
\newblock Tannaka duality and convolution for duoidal categories.
\newblock {\em Theory Appl. Categ.}, 28:No. 6, 166--205, 2013.

\bibitem{Bor:1994ii}
F.~Borceux.
\newblock {\em Handbook of categorical algebra. \textup2}, volume~51 of {\em
  Encyclopedia Math. Appl.}
\newblock Cambridge Univ. Press, Cambridge, 1994.

\bibitem{BruVir2007}
A.~Brugui{\`e}res and A.~Virelizier.
\newblock Hopf monads.
\newblock {\em Adv. Math.}, 215(2):679--733, 2007.

\bibitem{CLS2010}
D.~Chikhladze, S.~Lack, and R.~Street.
\newblock Hopf monoidal comonads.
\newblock {\em Theory Appl. Categ.}, 24:No. 19, 554--563, 2010.

\bibitem{For2004}
S.~Forcey.
\newblock Enrichment over iterated monoidal categories.
\newblock {\em Algebr. Geom. Topol.}, 4:95--119 (electronic), 2004.

\bibitem{Fra1976}
G.~A. Fraser.
\newblock The tensor product of distributive lattices.
\newblock {\em Proc. Edinburgh Math. Soc. (2)}, 20(2):121--131, 1976.

\bibitem{G}
R.~Garner.
\newblock Understanding the small object argument.
\newblock {\em Appl. Categ. Structures}, 17(3):247--285, 2009.

\bibitem{GL2016}
R.~Garner and I.~L{\'o}pez~Franco.
\newblock Commutativity.
\newblock {\em J. Pure Appl. Algebra}, 220(5):1707--1751, 2016.

\bibitem{Gol1999}
J.~S. Golan.
\newblock {\em Semirings and their applications}.
\newblock Kluwer Academic Publishers, Dordrecht, 1999.
\newblock Updated and expanded version of {{\i}t The theory of semirings, with
  applications to mathematics and theoretical computer science} [Longman Sci.
  Tech., Harlow, 1992; MR1163371 (93b:16085)].

\bibitem{H}
C.~Hermida.
\newblock Representable multicategories.
\newblock {\em Adv. Math.}, 151(2):164--225, 2000.

\bibitem{J}
P.~T. Johnstone.
\newblock {\em Topos theory}.
\newblock Academic Press [Harcourt Brace Jovanovich Publishers], London, 1977.
\newblock London Mathematical Society Monographs, Vol. 10.

\bibitem{JS}
A.~Joyal and R.~Street.
\newblock Braided tensor categories.
\newblock {\em Adv. Math.}, 102(1):20--78, 1993.

\bibitem{JT}
A.~Joyal and M.~Tierney.
\newblock An extension of the {G}alois theory of {G}rothendieck.
\newblock {\em Mem. Amer. Math. Soc.}, 51(309):vii+71, 1984.

\bibitem{Kel74}
G.~M. Kelly.
\newblock Doctrinal adjunction.
\newblock In {\em Category Seminar (Proc. Sem., Sydney, 1972/1973)}, pages
  257--280. Lecture Notes in Math., Vol. 420. Springer, Berlin, 1974.

\bibitem{K70}
A.~Kock.
\newblock Monads on symmetric monoidal closed categories.
\newblock {\em Arch. Math. (Basel)}, 21:1--10, 1970.

\bibitem{K71a}
A.~Kock.
\newblock Bilinearity and {C}artesian closed monads.
\newblock {\em Math. Scand.}, 29:161--174 (1972), 1971.

\bibitem{K71b}
A.~Kock.
\newblock Closed categories generated by commutative monads.
\newblock {\em J. Austral. Math. Soc.}, 12:405--424, 1971.

\bibitem{K72}
A.~Kock.
\newblock Strong functors and monoidal monads.
\newblock {\em Arch. Math. (Basel)}, 23:113--120, 1972.

\bibitem{LS2014}
S.~Lack and R.~Street.
\newblock A skew-duoidal {E}ckmann-{H}ilton argument and quantum categories.
\newblock {\em Appl. Categ. Structures}, 22(5-6):789--803, 2014.

\bibitem{Le}
T.~Leinster, editor.
\newblock {\em Higher operads, higher categories}, volume 298 of {\em London
  Mathematical Society Lecture Note Series}.
\newblock Cambridge University Press, Cambridge, 2004.

\bibitem{Li}
F.~E.~J. Linton.
\newblock Coequalizers in categories of algebras.
\newblock In {\em Sem. on {T}riples and {C}ategorical {H}omology {T}heory
  ({ETH}, {Z}\"urich, 1966/67)}, pages 75--90. Springer, Berlin, 1969.

\bibitem{McL}
S.~Mac~Lane.
\newblock {\em Categories for the working mathematician}.
\newblock Springer-Verlag, New York, second edition, 1998.

\bibitem{McC}
P.~McCrudden.
\newblock Opmonoidal monads.
\newblock {\em Theory Appl. Categ.}, 10:No. 19, 469--485, 2002.

\bibitem{MesWis2011}
B.~Mesablishvili and R.~Wisbauer.
\newblock Bimonads and {H}opf monads on categories.
\newblock {\em J. K-Theory}, 7(2):349--388, 2011.

\bibitem{M}
I.~Moerdijk.
\newblock Monads on tensor categories.
\newblock {\em J. Pure Appl. Algebra}, 168(2-3):189--208, 2002.
\newblock Category theory 1999 (Coimbra).

\bibitem{S}
G.~J. Seal.
\newblock Tensors, monads and actions.
\newblock {\em Theory Appl. Categ.}, 28:No. 15, 403--433, 2013.

\bibitem{St2012}
R.~Street.
\newblock Monoidal categories in, and linking, geometry and algebra.
\newblock {\em Bull. Belg. Math. Soc. Simon Stevin}, 19(5):769--821, 2012.

\bibitem{Szl}
K.~Szlach{\'a}nyi.
\newblock The monoidal {E}ilenberg-{M}oore construction and bialgebroids.
\newblock {\em J. Pure Appl. Algebra}, 182(2-3):287--315, 2003.

\end{thebibliography}
\end{document}